\documentclass[12pt]{article}
\usepackage{amsfonts}
\usepackage{pictex}

\def\baselinestretch{1} \topmargin -12pt \headsep 0pt \footskip 30pt
\def\medskipamount{12pt} \def\smallskipamount{6pt}

\newcounter{bitcount}
\newcommand{\bit}[1]{\addtocounter{bitcount}{1}\pagebreak[3]
\subsection{#1}\nopagebreak\setcounter{equation}{0}}

\renewcommand{\theequation}{\thesubsection .\arabic{equation}}
\renewcommand{\thesubsection}{\arabic{bitcount}}

\newcommand{\reb}[1]{\mbox{\bf (\ref{#1})}}
\newcommand{\re}[2]{\mbox{(\ref{#1}#2)}}

\catcode`\@=\active
\catcode`\@=11
\def\@eqnnum{\hbox to .01pt{}\rlap{\bf \hskip -\displaywidth(\theequation)}}
\catcode`\@=12

\begin{document}


\catcode`\@=\active
\catcode`\@=11
\newcommand{\nc}{\newcommand}


\nc{\bs}[1]{ \addvspace{\medskipamount} \pagebreak[3]
\refstepcounter{equation}
\noindent {\bf (\theequation) #1.} \begin{em} \nopagebreak}

\nc{\es}{\end{em} \par \addvspace{\medskipamount} } 

\nc{\ess}{\end{em} \par}

\nc{\br}[1]{ \addvspace{\medskipamount} \pagebreak[3]
\refstepcounter{equation} 
\noindent {\bf (\theequation) #1.} \nopagebreak}

\nc{\brs}[1]{ \pagebreak[3]
\refstepcounter{equation} 
\noindent {\bf (\theequation) #1.} \nopagebreak}

\nc{\er}{\par \addvspace{\medskipamount} }


\nc{\vars}[2]
{{\mathchoice{\mb{#1}}{\mb{#1}}{\mb{#2}}{\mb{#2}}}}
\nc{\C}{\mathbb C}
\nc{\Pj}{\mathbb P}
\nc{\Q}{\mathbb Q}
\nc{\R}{\mathbb R}
\nc{\Z}{\mathbb Z}
 

\nc{\oper}[1]{\mathop{\mathchoice{\mbox{\rm #1}}{\mbox{\rm #1}}
{\mbox{\rm \scriptsize #1}}{\mbox{\rm \tiny #1}}}\nolimits}
\nc{\coker}{\oper{coker}}
\nc{\ev}{\oper{ev}}
\nc{\Gr}{\oper{Gr}}
\nc{\Hom}{\oper{Hom}}
\nc{\LaGr}{\oper{LaGr}}
\nc{\NS}{\oper{NS}}
\nc{\OGr}{\oper{OGr}}
\nc{\Pic}{\oper{Pic}}
\nc{\Proj}{\oper{Proj}}
\nc{\Sec}{\oper{Sec}}
\nc{\Sym}{\oper{Sym}}

\nc{\operlim}[1]{\mathop{\mathchoice{\mbox{\rm #1}}{\mbox{\rm #1}}
{\mbox{\rm \scriptsize #1}}{\mbox{\rm \tiny #1}}}}

\nc{\bigtimes}{\operlim{$\times$}}

\nc{\PGL}[1]{{{\rm PGL(}#1{\rm )}}}
\nc{\SL}[1]{{{\rm SL(}#1{\rm )}}}
\nc{\U}[1]{{{\rm U(}#1{\rm )}}}

\nc{\Ga}{\Gamma}
\nc{\La}{\Lambda}
\nc{\la}{\lambda}
\nc{\si}{\sigma}

\nc{\bino}[2]{\mbox{\Large $#1 \choose #2$}}


\nc{\beqas}{\begin{eqnarray*}}
\nc{\ci}{{\cal I}}
\nc{\cH}{{\cal H}}
\nc{\cU}{{\cal U}}
\nc{\cC}{{\cal C}}
\nc{\co}{{\cal O}}
\nc{\cx}{{\C^{\times}}}
\nc{\down}{\Bigg\downarrow}
\nc{\downarg}[1]{{\phantom{}\Big\downarrow \raisebox{.4ex}{$\scriptstyle #1$}}}
\nc{\eeqas}{\end{eqnarray*}}
\nc{\emb}{\hookrightarrow}
\nc{\fp}{\mbox{     $\maltese$} \par \addvspace{\smallskipamount}}
\nc{\lra}{\longleftrightarrow}
\nc{\lrow}{\longrightarrow}
\nc{\mo}[1]{{\overline{M}_{0,#1}}}
\nc{\pf}{\noindent {\em Proof}}
\nc{\st}{\, | \,}
\nc{\uda}{\Big\updownarrow}
\nc{\udarg}[1]{\raisebox{.4ex}{$\phantom{\scriptstyle #1}$}
\Big\updownarrow\raisebox{.4ex}{$\scriptstyle #1$}}

\nc{\X}{{\cal X}}
\nc{\Y}{{\cal Y}}

\nc{\half}{\mathchoice{{\textstyle \frac{\scriptstyle 1}{\scriptstyle
2}}} {{\textstyle\frac{\scriptstyle 1}{\scriptstyle 2}}}
{\frac{\scriptscriptstyle 1}{\scriptscriptstyle 2}}
{\frac{\scriptscriptstyle 1}{\scriptscriptstyle 2}}}

\nc{\mod}{\, / \,}
\nc{\modd}{/}
\nc{\chow}{/ \! \! /}
\nc{\hilb}{/ \! \! / \! \! /}

\nc{\lp}{\raisebox{-.1ex}{\rm\large(}}
\nc{\rp}{\raisebox{-.1ex}{\rm\large)}}

\renewcommand{\labelenumi}{{\rm (\theenumi )}}

\hyphenation{col-line-a-tion col-line-a-tions para-met-riz-ing}

\catcode`\@=12



\noindent
{\LARGE \bf Complete collineations revisited}
\medskip \\ 
{\bf Michael Thaddeus } \smallskip \\ 
Department of Mathematics, Columbia University \\
2990 Broadway, New York, N.Y. 10027
\renewcommand{\thefootnote}{}
\footnotetext{Partially supported by NSF grants DMS--9500964 and
  DMS--9808529.} 

\setcounter{bitcount}{-1} 

\bit{Introduction}

This paper takes a new look at some old spaces.

The old spaces are the {\em moduli spaces of complete collineations},
introduced and explored by many of the leading lights of 19th-century
algebraic geometry, such as Chasles, Schubert, Hirst, and Giambelli.
They are roughly compactifications of the spaces of linear maps of a
fixed rank between two fixed vector spaces, in which the boundary
added is a divisor with normal crossings.  This renders them useful in
solving many enumerative problems on linear maps, and they are famous
as much for the intricacy of the resulting formulas as for the
elegance and symmetry of the underlying geometry.

The new look comes from some recent quotient constructions in
algebraic geometry: the {\em Chow quotient} and {\em inverse limit of
Mumford quotients}, introduced in 1991 by Kapranov, Sturmfels, and
Zelevinsky, and the {\em Hilbert quotient}, essentially due to
Bialynicki-Birula and Sommese in 1987.  These were motivated partly by
the search for a quotient more canonical than the Mumford or geometric
invariant theory quotient, which depends on the choice of a
linearization, and partly by the attractive polyhedral interpretations
possessed by all three quotients in the setting of toric geometry.
But for us, their chief interest arose later, in two papers of
Kapranov from 1993.  They showed, among many other things, that Chow
quotients, and several related operations, could be used to give
elegant new constructions of the moduli space $\mo{n}$ of stable
punctured curves of genus zero.

This paper will demonstrate that every one of Kapranov's constructions
has a counterpart for complete collineations.  In fact, one piece of
the puzzle was already in place: Vainsencher had already shown in 1984
that the complete collineations could be obtained by a blow-up
construction very similar to one given by Kapranov for $\mo{n}$.  But
all the other pieces also fit neatly.  For example, $\mo{n}$ is
constructed by Kapranov as the Chow quotient of a Grassmannian by a
$(\cx)^n$-action; likewise, the complete collineations are constructed
here as the Chow quotient of a Grassmannian by a $\cx$-action.

\br{Outline of the paper} 
Section \ref{review} of the paper reviews the quotient constructions
mentioned above, and Kapranov's work on $\mo{n}$.  Section
\ref{collineations} then reviews Vainsencher's work and states our
main theorem, which closely parallels Kapranov's.

Sections \ref{grass} and \ref{secant} are an apparent digression on
the varieties obtained as Mumford quotients of the Grassmannian by the
aforementioned $\cx$-action.  It turns out that they form a sequence
of smooth varieties obtained from a projective space by blowing up and
down the secant varieties of the Segre embedding, much like the
author's previous work \cite{thad-pair,thad-toric,thad-flip}.  This
identifies the inverse limit of Mumford quotients with Vainsencher's
construction, since his blow-up loci are the proper transforms of the
secant varieties.

Sections \ref{limit} through \ref{graphs} proceed with the proof of
the main theorem.  One key step, described in \S\ref{gelfand}, is a
version of the Gel'fand-MacPherson correspondence relating
quotients by a reductive group to quotients by a torus.  Another is
the construction of a morphism from the Chow quotient to any Mumford
quotient; this is due to Kapranov but re-proved in \S\ref{morphism}.

Section \ref{gw} is a brief remark explaining how the complete
collineations can also be interpreted as a space of stable maps to a
Grassmannian, and how this can be used to evaluate certain
Gromov-Witten invariants. 

The final three sections are appendices.  Sections \ref{quadrics} and
\ref{skew} deal with the special case when the two vector spaces in
the complete collineations are dual to each other.  In this case one
can restrict to the symmetric and anti-symmetric loci, to obtain
smaller spaces: the classical space of complete quadrics, and the
space of complete skew forms.  Everything stated for complete
collineations has a symmetric and anti-symmetric analogue in an
appealing way.  For example, the Grassmannian becomes a Lagrangian or
orthogonal Grassmannian, and the Segre embedding becomes a quadratic
Veronese or Pl\"ucker embedding.

Section \ref{symplectic}, although it discusses the original
motivation for the paper, has been banished to the end because it
breaks the flow.  It explains how the Chow quotient at the heart of
the paper may be understood in terms of a familiar object in
symplectic geometry, namely the space of broken Morse flows of the
moment map for a circle action.  Nothing from this section is needed
elsewhere in the paper, but it sheds some light on the geometric
intuition behind the Chow quotient.  \er

\br{Conventions} The base field is the field of complex numbers, or
any algebraically closed field of characteristic 0.  The
morphism from the Hilbert scheme to the Chow variety requires
characteristic 0, but many of the results not involving the Hilbert
scheme should hold in arbitrary characteristic. 

For a vector space $U$, $\Gr_i U$ denotes the Grassmannian of
$i$-dimensional subspaces of $U$.  The tautological subbundle and
quotient bundle over $\Gr_i U$ are denoted $S_U$ and $Q_U$,
respectively.  Also, $\Pj U$ denotes the space of lines in $U$, so
that $\Pj = \Gr_1$. 

The notation $X \hilb G$, $X \chow G$, $X \modd G$, for Hilbert, Chow
and Mumford quotients respectively, follows Kapranov \cite{kap-chow},
but it conflicts with many other papers where $X \chow G$ denotes a
Mumford quotient. 

The word {\em unstable} is used to mean ``not semistable''.

The term {\em complete collineation} is reserved by Laksov
\cite{bowdoin}, but not Vainsencher \cite{vain-clm}, for the case
where the two vector spaces are isomorphic.  Despite the loss of
precision, we prefer {\em collineation} for its euphony.  \er

\br{Acknowledgements} I am very grateful to the Institut
Mittag-Leffler for its warm hospitality in the spring of 1997; also to
William Fulton and Dan Laksov for their kind invitation to visit the
institute, and for helpful discussions on subjects related to this
paper. 

After this paper was written, I learned from M. M. Kapranov of the
previous work of Yu.\ A. Neretin.  This includes a result similar to the
identification of (\ref{clm}) and (\ref{gr}b) in the Main Theorem.
However, Neretin's work is not in algebraic geometry, but rather in the
category of metric spaces.

\bit{Review of algebraic quotient constructions \label{review}}

\brs{ Hilbert quotients}  Let $G$ be an algebraic group acting on a
projective variety $X$.  There is an open set $U \subset X$ such that
the orbit closures of points in $U$ form a flat family of subschemes
of $X$.  Indeed, let $\Delta \subset X \times X$ be the diagonal; the
closure $Y = \overline{G \cdot \Delta}$ is a projective variety, with
a morphism $Y \to X$ given by projection on the first factor.  For
dimension reasons, the fiber of this morphism over a generic $x \in X$
is simply the orbit closure $\overline{G \cdot x}$, an integral
subscheme of $X$.  Since $X$ is integral, there is a nonempty open set
$U \subset X$ such that $Y |_U$ is flat \cite[\S8]{cas}.

Hence there is a morphism $U \to \cH_X$, where $\cH_X$ is the
Hilbert scheme parametrizing subschemes of $X$.  Define the {\em
Hilbert quotient} $X \hilb G$ as the closure of the image of $U$ in
$\cH_X$.  This is a projective variety, and is independent of $U$.
It parametrizes a family of orbit closures, together with their flat
limits as schemes.

\br{Chow quotients} The family of schemes over the open set $U$ above
can also be viewed as a family of algebraic cycles.  Since $U$ is
reduced, this determines a morphism from $U$ to the Chow variety
parametrizing algebraic cycles on $X$ of the appropriate intersection
class \cite[I 3.10, 3.21]{kollar}.  Define the {\em Chow quotient} $X
\chow G$ as the closure of the image of $U$ in the Chow variety.
Again, this is a projective variety, and is independent of $U$.  It
parametrizes a family of orbit closures, together with their limits as
cycles.

In fact, Mumford \cite[\S5.4]{mfk} showed that in characteristic 0,
there is a morphism from the induced reduced subscheme of every
component of the Hilbert scheme to the appropriate Chow variety.
Hence there is a canonical morphism from the Hilbert quotient to the
Chow quotient.

\br{Inverse limits of Mumford quotients\label{inv-lim}} Suppose now
that $X$ is normal and $G$ is reductive.  Geometric invariant theory
associates a {\em Mumford quotient} $X \modd G \, (L)$ to every ample
(or even base-point free) linearization $L$ of the $G$-action.  Here a
{\em linearization} is a line bundle on $X$ endowed with a linear
action of $G$ lifting the action on $X$.  The quotient is defined to
be $\Proj \bigoplus_k H^0(L^k)^G$, where the superscript denotes the
$G$-invariant part.  It is the categorical quotient of an open subset
of $X$, called the {\em semistable set}.

The Mumford quotients are familiar, but depend on the non-canonical
choice of a linearization.  Less familiar, but more canonical, is
their inverse limit, which will be our focus of interest. 

The Mumford quotients are given the structure of an inverse system in
the following way.  First, let the $G$-{\em N\'eron-Severi group}
$\NS^G(X)$ be the group of $G$-algebraic equivalence classes of
linearizations on $X$, where $G$-{\em algebraic equivalence} is the
obvious generalization of ordinary algebraic equivalence to
linearizations.  Then there is a short exact sequence
$$0 \lrow \chi(G) \otimes \Q \lrow \NS^G(X) \otimes \Q \lrow \NS(X)
\otimes \Q \lrow 0,$$
where $\chi(G) = \Hom(G, \cx)$.  Moreover, up to
isomorphism the quotient $X \modd G\,(L)$ depends only on the class of
$L$ in $\Pj (\NS^G(X) \otimes \Q)$.  In a slight abuse of terminology,
an element of this rational projective space will be called a {\em
fractional linearization}.  Second, define a {\em chamber structure}
on a domain in a rational projective space to be a partition that is
locally a linear projection of the faces of a rational polytope.  Then
there is a chamber structure on the set of base-point free
linearizations in $\Pj (\NS^G(X) \otimes \Q)$, having only finitely
many chambers, such that up to isomorphism, the quotient $X \modd G
(L)$ depends only on the chamber containing $L$.  Third, taking
closures gives the set of chambers the structure of a directed set.
For any arrow in this directed set, there is a natural morphism of the
corresponding quotients.  This makes the collection of all nonempty
Mumford quotients $X \modd G\,(L)$ into an inverse system.  Because the
chambers are finite in number, the inverse limit is a projective
scheme of finite type over $\C$.

The statements in the previous paragraph are proved in the works of
Dolgachev-Hu \cite{dh} and the author \cite{thad-flip}; in the torus
case, many of them are due to Brion-Procesi \cite{bp}.  However, only
one of them is substantially difficult to prove, namely the finiteness
of the chambers.  Everything else in some sense follows from the
Hilbert-Mumford numerical criterion \cite{mfk}.  Moreover, in the
cases treated herein, the finiteness is also fairly easy.  So the
inverse limits we need can be defined without invoking the general
theory.  \er

The Hilbert quotient, the Chow quotient, and the limit of the Mumford
quotients are all canonically associated to a normal variety with a
reductive group action.  In the case of torus actions on toric
varieties, they have been studied by Kapranov-Sturmfels-Zelevinsky
\cite{ksz}, where they are generally all different.  In particular,
the Chow and Hilbert quotients are irreducible by definition, while
the inverse limit may be reducible even in a very simple example
\cite[1.11]{thad-flip}.  However, in any case the limit of Mumford
quotients has a {\em distinguished component} containing the quotients
of that open set of points which are semistable whenever the quotient
is nonempty.
 
\medskip
 
Kapranov \cite{kap-chow,kap-vero} went on to make a detailed study of
these quotients in a specific case.  We conclude this section by
stating his results, which were a principal inspiration for the
present work.  They give several remarkable constructions of the
moduli space of stable curves of genus 0.

For any $n$, let $\Ga_n$ be a set of $n$ points in general position in
$\Pj^{n-2}$.  Note that all such sets are projectively equivalent.
Define varieties $X_{n,k}$ recursively as follows: $X_{n,1} =
\Pj^{n-3}$, and $X_{n,k+1}$ is the blow-up of $X_{n,k}$ along the
proper transforms of all the $k-1$-planes spanned by $k$ points in
$\Ga_{n-1}$.  Kapranov \cite{kap-chow} is very careful about the order
in which these $k$-planes are blown up.  However, as D. Thurston has
pointed out, when two smooth subvarieties are transverse, the order in
which their proper transforms are blown up does not matter.  All the
$k$-planes in question are pairwise transverse, so blowing them up in
any order yields the same variety.

\bs{Theorem (Kapranov)} The following are isomorphic:

\begin{enumerate}

\item The moduli space $\mo{n}$ of stable $n$-punctured curves of
genus 0;

\item The variety $X_{n,n-2}$ defined above;

\item The closure {\rm (a)} in the Hilbert scheme, or {\rm (b)} in the
Chow variety of $\Pj^{n-2}$, of the locus of rational normal curves
containing $\Ga_n$;

\item The {\rm (a)} Chow quotient, or {\rm (b)} inverse limit of
Mumford quotients $(\Pj^1)^n \modd \SL{2}$; 

\item The {\rm (a)} Hilbert quotient, or {\rm (b)} Chow quotient, or
  {\rm (c)} inverse limit of Mumford quotients $\Gr_2 \C^n \mod
  (\cx)^n$, where $(\cx)^n$ acts in the natural way. \fp
\end{enumerate} 
\es

\bit{Complete collineations \label{collineations}}
\brs{Definitions\label{defns}} 
The spaces of complete collineations were introduced and extensively
studied at the end of the nineteenth century.  The literature on them
is vast and we content ourselves with a very sketchy review.  A
comprehensive history, both ancient and modern, is recounted by Laksov
\cite{nice,bowdoin}.

Let $V$, $W$ be finite-dimensional vector spaces, and let $u \leq
\min(\dim V, \dim W)$.  Over $\Gr_u V$, there is a rational map
$$\Pj \Hom(S_V,W) \dasharrow \bigtimes_{i=0}^u \Pj \Hom (\La^i S_V,
\La^i W)$$
given by $f \mapsto (\La^i f)$.  The {\em space of rank $u$
complete collineations} is defined to be the closure of the image.
By definition, this compactifies the space of isomorphisms between
$u$-dimensional subspaces of $V$ and $W$, up to scalars.  When $u =
\dim V = \dim W$, for example, it is a compactification of $\PGL{u}$.

The name given to this space suggests that there should exist an object
called a rank $u$ complete collineation of which it is the moduli
space.  This is indeed true: a {\em rank $u$ complete collineation} $V
\to W$ is a $u$-dimensional subspace $U$ of $V$, together with a
finite sequence of nonzero linear maps $f_i$, where $f_1: U \to W$,
$f_{i+1} : \ker f_i \to \coker f_i$, and the last $f_i$ has maximal
rank.

However, the notion of families of such objects is rather delicate; it
is meticulously treated by Kleiman-Thorup \cite{kt} and Laksov
\cite{lak-arkiv}.  Roughly speaking, where $f_i$ drops rank, its
derivative should be $f_{i+1}$.  Rather than struggle with this
interpretation, we will make contact with the complete collineations
via an explicit construction of the moduli space, due to Vainsencher
in 1984 \cite{vain-clm}. \er

\br{Vainsencher's construction\label{vc}}
Given $V$, $W$, and $u$ as above, define a sequence of varieties
$\X^u_k(V,W)$, or just $\X_k$ for short, by the following recursion.
Let $\X_1$ be the projective bundle $\Pj \Hom(S_V,W)$ mentioned above,
and let $\X_k$ be the blow-up of $\X_{k-1}$ at the proper transform of
the locus of linear maps $S_V \to W$ of rank $k-1$.  Vainsencher
proved that $\X^u_u(V,W)$ is isomorphic to the space of rank $u$
complete collineations.

If $u= \dim V$, for example, then $\X_1$ is just the projective space
$\Pj \Hom(V,W)$, and the blow-up loci are the proper transforms of the
secant varieties of the Segre embedding $\Pj V^* \times \Pj W \emb \Pj
\Hom(V,W)$.  \er

The similarity of this recursive construction to Kapranov's is
striking.  Our main result shows that this analogy goes much
further. 

\bs{Main Theorem} The following are isomorphic:

\begin{enumerate}

\item The moduli space of rank $u$ complete collineations $V \to W$;
\label{clm}

\item The variety $\X^u_u(V,W)$ defined above;
\label{vain}

\item The closure {\rm (a)} in the Hilbert scheme, or {\rm (b)} in the
Chow variety of $\Pj V \times \Pj W$, of the locus of graphs of
linear injections $U \to W$, where $U \subset V$ ranges over all
$u$-dimensional subspaces;
\label{graph}

\item The {\rm (a)} Chow quotient, or {\rm (b)} inverse limit of
Mumford quotients $\Pj \Hom(U, V^*) \times \Pj \Hom (U,W) \mod
\SL{U}$, where $U$ is a fixed vector space of dimension $u$;
\label{proj}

\item The {\rm (a)} Hilbert quotient, or {\rm (b)} Chow quotient, or
  {\rm (c)} inverse limit of Mumford quotients $\Gr_u (V \oplus W)
  \modd \cx$, where $\cx$ acts with weight $1$ on $V$ and $-1$ on $W$.
\label{gr}

\end{enumerate}

\es

As we have seen, the equivalence of \re{clm}{} and \re{vain}{} is
Vainsencher's, but the rest appears to be new.  The proof of the
remaining statements occupies sections \ref{grass} through
\ref{graphs}.  The diagram below shows which equivalences are proved
in which sections.

$$\begin{array}{ccccccccccc}
 & & & & & &\re{gr}{a}& & & & \\
 & & & & & &\udarg{\ref{gr-chow/hilb}} & & & & \\
\re{clm}{} & \stackrel{\ref{vc}}{\lra} & \re{vain}{} & 
\stackrel{\ref{vain/gr-mum}}{\lra} & \re{gr}{c} &
\stackrel{\ref{gr-chow/mum}}{\lra} 
& \re{gr}{b} 
& \stackrel{\ref{gr/graph}}{\lra} & \re{graph}{b} 
& \stackrel{\ref{graph-chow/hilb}}{\lra} & \re{graph}{a} \\
 & & & & \udarg{\ref{mum-gr/proj}} 
 & & \udarg{\ref{chow-gr/proj}} & & & & \\
 & & & & \re{proj}{b} & & \re{proj}{a} & & & &
\end{array}
$$

\bit{Quotients of the Grassmannian by the multiplicative group
  \label{grass}}
As in \S\ref{collineations} above, $V, W$ are finite-dimensional
vector spaces, and $u$ is a positive integer such that $u \leq
\min(\dim V, \dim W)$.  Let $T = \cx$ act on $V \oplus W$ by $\la
(v,w) = (\la v, \la^{-1} w)$; this induces a $T$-action on $Y = \Gr_u
(V \oplus W)$.  This section concerns the Mumford quotients of $Y$ by
this action.  In fact, all the constructions are sufficiently natural
that $V$ and $W$ could, more generally, be vector bundles over a
common base.

The $T$-action is canonically linearized on the ample line bundle
$\det Q_{V\oplus W}$.  However, the canonical linearization is not the
only one, because the group of characters $\chi(T)$ is not trivial,
but isomorphic to $\Z$.  For any $t \in \chi(T) \otimes \Q \cong \Q$,
the canonical linearization $L$ can be tensored by the fractional
character $t$ to produce a new fractional linearization $L(t)$.
Indeed, by the short exact sequence of \reb{inv-lim}, any ample
linearization is a positive tensor power of some $L(t)$.

For $\si \in \Q$, say $U \in Y$ is {\em $\si$-semistable} if it is
semistable for the linearization $L(u - 2\si)$, and let $Y_\si$ be the
corresponding Mumford quotient.  We aim to compare the quotients for
different $\si$, using the results of \cite{hu, thad-flip}.  Actually,
in the torus case those results are only slight refinements of the
work of Brion-Procesi \cite{bp}.

A straightforward application of the Hilbert-Mumford numerical
criterion \cite{mfk} yields the following result.

\bs{Proposition\label{stability}}
A point $U \in Y$ is 
{\renewcommand{\labelenumi}{{\rm (\alph{enumi})}}
\begin{enumerate}
\item $\si$-semistable if and only if $\dim U \cap V
\leq u - \si$ and $\dim U \cap W \leq \si$;
\item $\si$-stable if and
only if both inequalities are strict;
\item identified in the quotient with
$(U \cap V) \oplus (U + V)/V$ if equality holds in the first case, and
with $(U + W)/W \oplus (U \cap W)$ if it holds in the second. \fp
\end{enumerate}}
\es

\bs{Corollary}
{\renewcommand{\labelenumi}{{\rm (\alph{enumi})}}
\begin{enumerate}
\item For $\si \not\in [0,u]$, $Y_\si = \emptyset$;

\item For $\si \in (0,1)$, $Y_\si$ is the
projective bundle $\Pj \Hom(S_V,W)$ over $\Gr_u V$;

\item For $\si \in (u-1,u)$, $Y_\si$ is the
projective bundle $\Pj \Hom(S_W,V)$ over $\Gr_u W$. \fp
\end{enumerate}}
\es

The semistability condition, and hence the quotient $Y_\si$, is
locally constant for $\si \not\in \Z$.  At an integer value, on the
other hand, the semistability condition changes.  The chamber
decomposition of \reb{inv-lim} therefore consists simply of
the integer points in $[0,u]$ and the intervals between them.  The
quotients for $\si \in [1, u-1]$ can then be described using Theorem
4.8 of \cite{thad-flip}, which explains how the quotient changes when
a wall is crossed.  Its hypotheses are readily verified, and it
directly implies the theorem below.

Fix an integer $k \in [1, u-1]$.  Write $Y^-_k$ for $Y_{k-\half}$ and
$Y^+_k$ for $Y_{k+\half}$.  Also let $Z^0_k$ be the locus in
$Y_k$ consisting of points unstable for $\si \neq k$, let
$Z^-_k$ be the locus in
$Y^-_k$ consisting of points unstable for $\si > k$, and let 
$Z^+_k$ be the locus in
$Y^+_k$ consisting of points unstable for $\si < k$.

\bs{Theorem\label{up&down}}
With the above notation, 
{\renewcommand{\labelenumi}{{\rm (\alph{enumi})}}
\begin{enumerate}
\item $Z^\pm_k$ are naturally projective bundles over
$Z^0_k$: indeed there are natural isomorphisms 
\begin{eqnarray*}
Z^0_k & = & \Gr_{u-k} V \times \Gr_k W, \\ 
Z^-_k & = & \Pj \Hom (S_W,Q_V), \\ 
Z^+_k & = & \Pj \Hom (S_V,Q_W);
\end{eqnarray*}
\item the birational morphisms of quotients $Y^-_k \to Y_k$ and
$Y^+_k \to Y_k$ induced by the inclusions $Y^{ss}(k \pm \half)
\subset Y^{ss}(k)$ restrict to the projections $Z^\pm_k \to Z^0_k$,
and are isomorphisms elsewhere;
\item the fibered product $Y^-_k \times_{Y_k} Y^+_k$ is
naturally isomorphic to the blow-up of $Y^-_k$ along
$Z^-_k$, and to the blow-up of $Y^+_k$ along
$Z^+_k$. \fp
\end{enumerate}}
\es

There is therefore a diagram of morphisms, all birational except in
the bottom corners:
{\def\arraycolsep{2pt} 
$$\begin{array}{ccccccccccccccccc}
& & & &  \tilde{Y}_1 & & & & \tilde{Y}_2 & & & &  \tilde{Y}_3 & & & & \\
& & & \swarrow & & \searrow & & \swarrow & & \searrow & & \swarrow & &  \searrow & & & \\
& & Y_{\half} & & & & Y_{1\half} & & & & Y_{2\half} 
& & & & \, \cdots \, & & \\
& \swarrow & & \searrow  & & \swarrow & & \searrow & & \swarrow & & \searrow 
& & \swarrow & & \searrow & \\
\Gr_u V & & & & Y_1 & & & & Y_2 & & & & Y_3 & & & &\Gr_u W;
\end{array}$$}here $\tilde{Y}_i$ denotes the fibered product, 
and the morphisms in the top row are the blow-ups of $Z^\pm_k$.

\bit{Secant variety interpretation\label{secant}} 
The results of the previous section have an appealing geometric
interpretation.  Denote by $\Sec_k(S_V,W)$ the locus in $Y_{\half} =
\Pj \Hom (S_V,W)$ consisting of linear maps of rank $\leq k$.  This is
precisely the variety of secant $(k-1)$-planes to the relative Segre
embedding of $\Pj S_V^* \times \Pj W$ in $\Pj(S_V^* \otimes \Pj W)$.

\bs{Proposition\label{blow-up/secant}}
The blow-up locus $Z^-_k$ is the proper transform of $\Sec_k(S_V,W)$.
\es

\pf.  By definition $Z^-_k$ consists of (quotients of) points unstable
for all $\si > k$.  By Proposition \ref{stability} this is the
quotient of the locus of subspaces $U \in Y$ satisfying $\dim (U \cap
V) \geq u-k$.  On the other hand, $Y^{ss}(\half)$ is simply the locus
of graphs of linear maps from a $u$-dimensional subspace of $V$ to
$W$.  These two loci meet in the locus of such graphs of rank $\leq
k$.  Its quotient $\Sec_k(S_W,V)$ is therefore the proper transform of
$Z^-_k$ in $Y_{\half}$.  \fp

A convenient example to keep in mind is when $\dim V = u$, so that
$\Gr_u V$ is a point.  Then $Y_{\half}$ is simply the projective space
$\Pj \Hom(V,W)$, and the locus of maps of rank $\leq i$ is the variety
of secant $(i-1)$-planes to the image of the Segre embedding $\Pj V^*
\times \Pj W \emb \Pj \Hom(V,W)$.  Theorem \ref{up&down} therefore
shows that the proper transforms of these secant varieties can be
blown up and down in turn, to produce a sequence of smooth varieties.
When $\dim W = \dim V$ as well, then $Y_{u-\half}$ is also a
projective space, and the rational map $Y_{\half} \dasharrow
Y_{u-\half}$ is the Cremona transformation given by inversion of the
linear map; the diagram factors it into a sequence of blow-ups and
blow-downs.

There are several other cases where the secant varieties of some
variety $M$ in projective space may be blown up and down in this way.
For example, according to the author's previous work \cite{thad-pair},
it is possible when $M$ is a smooth curve of genus $g$ embedded by a
complete linear system of degree $> 2g-2$.  It is also possible when
$M$ consists of $< n+1$ \cite{kap-chow,thad-toric} or $n+2$
\cite[\S6]{thad-flip} general points in $\Pj^n$.  And we will see in
\S\S\ref{quadrics} and \ref{skew} that it is possible for the
quadratic Veronese embedding of $\Pj^n$ and for the Pl\"ucker
embedding of $\Gr_2 \C^n$, just by taking appropriate cross-sections
of the Segre embedding.

It is quite an interesting and challenging problem to understand how
generally this phenomenon occurs.  On the one hand, the list assembled
above looks promising, but on the other hand, simple examples (such as
that of 6 general points in $\Pj^3$) show that the varieties involved
may not be projective.  In that case, the construction of blow-downs
becomes extremely problematic.  We hope to return to this fascinating
subject in the future.

\bit{The inverse limit of Mumford quotients of the
  Grassmannian\label{limit}} 

The Main Theorem is concerned not with any single Mumford quotient,
but with their inverse limit.  In the case of the action of $T = \cx$
on $Y = \Gr_u(V \oplus W)$ from \S\ref{grass}, this can be constructed
as follows.

For any positive integer $k < u$, define a sequence of varieties
$\Y^u_k(V,W)$, or just $\Y-k$ for short, and morphisms $\phi_k: \Y_k
\to Y_k$ by the following recursion.  Let $\Y_1$ be $Y_{\half}$,
and given $\Y_k$, let $\Y_{k+1}$ be the fibered product
$$\Y_k \times_{Y_k} Y_{k+\half}.$$
Then $\Y^u_u(V,W)$ is precisely the inverse limit of all the Mumford
quotients, and will shortly be shown isomorphic to the space of complete
collineations. 

The intricate geometry of these spaces arises from a recursive
structure: each one contains similar spaces of lower
dimension.  

\bs{Proposition\label{recursion}} 
Over $Z^0_k = \Gr_{u-k} V \times \Gr_k W$, the restriction of $\phi_k$
is the projection from $\Y^k_k(Q_V,S_W)$, that is, the locally
trivial family whose fiber over $(V',W')$ is $\Y^k_k(V/V', \, W')$.  
\es

\pf.  First let us determine $\phi_k^{-1}(V',W')$ for any single
element $(V',W') \in Z^0_k$.  Let
$$Y' = \{ U \in Y \st V' \subset U \subset V \oplus W' \}.$$
This is a subvariety of $Y$ naturally isomorphic to $\Gr_k(V/V' \oplus
W')$.  Proposition \ref{stability} implies that the quotient $Y'_\si
\subset Y_\si$ is nonempty if and only if $\si \in [0,k]$, equals the
single point $(V',W')$ if $\si = k$, and is preserved by taking the
image in the natural maps $Y^\pm_k \to Y_k$, as well as the {\em
inverse} image in one direction $Y^-_k \to Y_k$ only.  It
immediately follows that $\phi_k^{-1}(V',W') = \Y^k_k(V/V', W')$.

To show that the collection of these spaces forms the locally trivial
family in the statement is straightforward: just consider, in place of
$Y'$, the locally trivial family of Grassmannians $\Gr_k(Q_V \oplus
S_W)$ parametrized by $Z^0_k$ and the natural map from this family to
$Y$.  \fp

In the Morse-theoretic terms outlined in \S\ref{symplectic}, $Y'$
is the closure of the upward Morse flow from $(V',W')$.

\br{Proof that \re{gr}{c} = \re{vain}{}\label{vain/gr-mum}}
We now show that the inverse limit $\Y^u_u(V,W)$ constructed 
above is isomorphic to Vainsencher's blow-up $\X^u_u(V,W)$.

By definition $\X_1 \cong \Y_1$.  By induction on $k$, it suffices to
assume $\X_k \cong \Y_k$ and prove $\X_{k+1} \cong \Y_{k+1}$.  Now
$\X_{k+1}$ is the blow-up of $\X_k$ along the proper transform of the
locus of rank $k$ linear maps, which by Proposition
\ref{blow-up/secant} is the proper transform of $Z^-_k$ in $\X_k$.  On
the other hand, by Theorem \ref{up&down} the fibered product $Y_{k-\half}
\times_{Y_k} Y_{k+\half}$ is the blow-up of $Y_{k-\half}$ at $Z^-_k$.
Taking the fibered product with $\Y_{k-1}$ over $Y_{k-1}$ shows that
$\Y_{k+1}$ is the blow-up of $\Y_k$ at $\Y_{k-1} \times_{Y_{k-1}}
Z^-_k$, which is the total transform of $Z^-_k$ in $\Y_k$.  It
therefore suffices to show that the total transform equals the proper
transform, or equivalently, that the total transform is an integral
scheme.

By Theorem \ref{up&down}, $Z^-_k$ is the inverse image of $Z^0_k$ in
the natural map $Y_{k-\half} \to Y_k$.  Its total transform is
therefore nothing but $\phi_k^{-1} (Z^0_k)$.  But this was shown in
Proposition \ref{recursion} to be a locally trivial fibration over
$\Gr_{u-k} V \times \Gr_k W$ with fiber $\Y^k_k(V/V', \, W')$.  By
induction on $u$, this may be assumed to be integral, and that completes
the proof.  \fp \er

\bit{The Gel'fand-MacPherson correspondence\label{gelfand}}

To relate the quotients of the Grassmannian to those of $\Pj \Hom
(U,V) \times \Pj \Hom (U,W)$, we use a version of the so-called
{\em Gel'fand-MacPherson correspondence} \cite{gm}.  The argument is also
reminiscent of the author's trick \cite[3.1]{thad-flip} for converting
a quotient by a reductive group to one by a torus.  It is carried out
here first for Mumford quotients, then for Chow quotients.

\br{Proof that \re{gr}{c} = \re{proj}{b}\label{mum-gr/proj}} 
Let $T = \cx$ act on $Y = \Gr_u(V \oplus W)$ as before.  We wish to
show that the inverse limit of Mumford quotients $Y \modd T$ is
isomorphic to that of $\Pj \Hom (U,V) \times \Pj \Hom (U,W) \modd
\SL{U}$.  It suffices to show that the inverse systems of Mumford
quotients are isomorphic.  We will see that both are isomorphic to the
inverse system of Mumford quotients $\Pj \Hom(U, V \oplus W) \modd \lp
T \times \SL{U}\rp$.

There is an unique ample linearization of
the $\SL{U}$-action on $\Pj \Hom(U, V \oplus W)$, and the Mumford
quotient is the Grassmannian $Y$.  It is an easy
exercise to check that the descent on linearizations induces an
isomorphism 
$$\NS_\Q^T(Y) \cong \NS_\Q^{T \times \SL{U}}\lp\Pj
\Hom(U, V \oplus W)\rp.$$  
It follows from the definition of Mumford quotients that 
$$\Pj \Hom(U, V \oplus W) \modd \lp T \times \SL{U}\rp \cong Y \modd
T$$
naturally, if the linearization on $Y$ is taken to be the descent
of that on the projective space.  Indeed, the isomorphism is induced
by a morphism of the corresponding semistable sets.  The morphisms in
the inverse systems are induced by inclusions of semistable sets, so
they commute with these isomorphisms.  Hence the inverse systems of $Y
\modd T$ and of \linebreak[4]
$\Pj \Hom(U, V \oplus W) \modd \lp T \times \SL{U}\rp$ are
isomorphic.

Similar arguments apply to the $T$-action on $\Pj \Hom(U, V \oplus
W)$.  The linearization is no longer unique, since it can be twisted
by a character of $T$.  But every nonempty quotient is naturally
isomorphic to $\Pj \Hom(U,V) \times \Pj \Hom(U,W)$.  (Actually, this
is not strictly true, as the special cases $\si = 0$ and $\si = u$
have just $\Pj \Hom(U,V)$ and $\Pj \Hom(U,W)$ as their respective
quotients.  However, it does no harm to pretend that the quotient
remains $\Pj \Hom(U,V) \times \Pj \Hom(U,W)$, since further dividing
by the $\SL{U}$-action gives $\Gr_uV$ and $\Gr_uW$ in either case.)
The descent again induces an isomorphism on the space of fractional
linearizations.  Hence the inverse systems of $\Pj \Hom(U, V \oplus W)
\modd \lp T \times \SL{U}\rp$ and of $\Pj \Hom(U,V) \times \Pj \Hom(U,W)
\modd \SL{U}$ are isomorphic.  \fp \er

\bs{Lemma\label{useful}}
Let $A$ be a cycle on $M \times N$ such that for any $x \in N$, the
naive intersection $A \cap \lp M \times \{ x \}\rp$ has the expected
dimension.  Then $A$ determines a morphism from $N$ to the appropriate
Chow variety of $M$.  
\es

\pf.  This is proved by Samuel \cite[II 5.12, 6.8]{samuel}, or, in
more modern language, by Koll\'ar \cite[I 3.10, 3.21]{kollar}.  \fp 

\br{Proof that \re{gr}{b} = \re{proj}{a}\label{chow-gr/proj}}
Consider the diagram 
{\def\arraycolsep{1pt}
$$\begin{array}{rcl}
\Pj \lp\co(1,0) \oplus \co(0,1)\rp & & \Pj \Hom(U,S_{V\oplus W}) \\
\down \phantom{xxxxxi..} & 
\begin{array}{ccc} 
\searrow & & \swarrow \\
 & \Pj\Hom (U,V \oplus W) & \\
\swarrow && \searrow 
\end{array}
& \phantom{xxxxx} \down \\
\Pj \Hom (U,V) \times \Pj \Hom (U,W) & & \Gr_u (V \oplus W).
\end{array}$$}

Here $\Pj \lp\co(1,0) \oplus \co(0,1)\rp$ refers to the $\Pj^1$-bundle
over $\Pj \Hom (U,V) \times \Pj \Hom (U,W)$, and the diagonal arrows
in the second row are only rational maps, geometric quotients of open
subsets by $T$ and $\SL{U}$, respectively.  

Certainly the universal cycle on $(Y \chow T) \times Y$ satisfies the
lemma.  It may be pulled back by the flat morphism in the right-hand
column.  The pulled-back cycle on $(Y \chow T) \times \Pj
\Hom(U,S_{V\oplus W})$ again satisfies the lemma, and over a general
point in $Y \chow T$, its fiber is an $\SL{U} \times T$-orbit
closure.  This cycle can be pushed forward to $(Y \chow T) \times
\Pj \Hom (U, V \oplus W)$.  The pushed-forward cycle satisfies the
lemma, since a birational morphism does not increase the codimension
of any subvariety.  But it is also $T$-invariant, so it descends to
a cycle on the geometric quotient $Y \chow T \times \Pj \Hom (U,V)
\times \Pj \Hom (U,W)$.  This again satisfies the lemma, hence
determines a morphism from $Y \chow T$ to a Chow variety of $\Pj \Hom
(U,V) \times \Pj \Hom (U,W)$.  Over a general point in $Y \chow T$,
its fiber is an $\SL{U}$-orbit closure.  Its image is therefore the
Chow quotient.

A similar argument using the $\Pj^1$-bundle in the diagram produces a
morphism in the opposite direction, which is clearly inverse to the
first one generically.  Since both spaces are integral, it must be the
inverse everywhere.  \fp \er

\bit{The morphism from Chow to Mumford quotients \label{morphism}}

The following theorem is essentially that proved by Kapranov
\cite[0.4.3]{kap-chow}, but with a few additions and a
different proof.

\bs{Theorem\label{chow->mum}}
Let $X$ be a projective variety, $G$ a reductive group acting on $X$.
For any linearization $L$ with $X^{ss}(L) \neq
\emptyset$, and any cycle $Y$ appearing in the Chow quotient $X \chow
G$:
{\renewcommand{\labelenumi}{{\rm (\alph{enumi})}}
\begin{enumerate}
\item $Y \cap X^{ss}(L) \neq \emptyset$ and is contained in a unique
$L$-semistable equivalence class of orbits;  
\item assigning this class to
$Y$ defines a morphism $X \chow G \to X \modd G\,(L)$, which is
birational if $X^{s}(L) \neq \emptyset$;  
\item if $L^+$ and
$L^0$ are two linearizations such that $X^{ss}(\!+\!) \subset X^{ss}(0)$,
then the diagram below commutes.
\end{enumerate}}
$$
\begin{array}{rcl}
& X \chow G & \\
\swarrow & & \searrow \\
X \modd G \, (\! + \! ) & \lrow & X \modd G \, (0)
\end{array}
$$
\es

\pf.  Denote $\cC$ the Chow variety containing the orbit closure of a
generic point in $X$.  The universal cycle $\cU$ on $C \times
X$ restricts via Chow pull-back to a universal cycle over $(X \chow G)
\times X$.  It is proper, $G$-invariant, irreducible, has multiplicity
1, and surjects onto both factors.  It therefore descends to a cycle
$D$ on $(X \chow G) \times (X^{ss} \modd G)$.  Denote $E$ the open set in
$X$ where the orbit closure has the generic homology class.  There is
a natural morphism $E \to X \chow G$, but also $E \cap X^{ss} \neq
\emptyset$ since both are nonempty and open.

For any $x \in X \chow G$ in the image of $E \cap X^{ss}$, $D \cap \lp
\{ x \} \times X^{ss} \modd G \rp$ clearly consists of a single point.
Moreover, for any $x \in X \chow G$ whatsoever, $D \cap \lp \{ x \}
\times X^{ss} \modd G\rp$ is finite: it consists exactly of the
semistable equivalence classes of irreducible components of $\cU \cap
\lp \{ x \} \times X^{ss}\rp$.  Lemma \ref{useful} then implies that
$D$ is the Chow pull-back of the universal cycle by a morphism from $X
\chow G$ to the Chow variety of 0-dimensional, length 1 cycles in
$X^{ss} \modd G$, which is of course $X^{ss} \modd G$ itself.  So $D
\cap \lp \{ x \} \times X^{ss} \modd G \rp$ consists of exactly one
point for any $x \in X \chow G$ at all.  There is hence a unique
semistable equivalence class in each cycle of the Chow quotient.
Furthermore, the morphism $X \chow G \to X^{ss} \modd G$ is clearly an
isomorphism over the open set $E \cap X^s \modd G$, so it is
birational if $X^s$ is nonempty.  This completes the proof of the
first two statements.

The third statement is easy after noting that, since the inclusion $(X
\chow G) \times X^{ss}(\!+\!) \subset (X \chow G) \times X^{ss}(0)$
induces an inclusion of universal cycles, the quotient cycle $D^+
\subset (X \chow G) \times X \modd G\,(+)$ is the inverse image of $D^0
\subset (X \chow G) \times X \modd G\,(0)$.  \fp

\bs{Corollary\label{corollary}}
For any reductive group action on a projective variety, there exists a
morphism from the Chow quotient to the inverse limit of Mumford
quotients. \fp
\es

At this point the reader may be thinking it likely that the Chow
quotient is always isomorphic to the inverse limit of Mumford
quotients.  In fact, this is far from being true.  A counterexample is
furnished by the action of $\SL{2}$ on $\Sym^n \Pj^1 = \Pj^n$.  The
Chow quotient parametrizes those cycles appearing in the quotient of
the universal cycle over $\mo{n} \times (\Pj^1)^n$ by the natural
action of the symmetric group $\Sigma_n$.  It is not hard to see that
the Chow quotient differs from both $\mo{n} \mod \Sigma_n$ and $\Pj^n
\modd \SL{2}$, and indeed is not normal at least when $n \geq 8$.
This may be the example alluded to in Remark 0.4.10 of Kapranov
\cite{kap-chow}.

\smallskip

However, let us return to the case at hand: $Y = \Gr_u(V \oplus W)$,
$T=\cx$.  It is now easy to describe those cycles belonging to the
Chow quotient $Y \chow T$.  The proof is straightforward using
Proposition \ref{stability} and the fact, stated in Theorem
\ref{chow->mum}, that any cycle in the Chow quotient intersects
$X^{ss}(\si)$ in exactly one $\si$-semistable equivalence class,
whether $\si$ is integral or not.

\bs{Proposition\label{nodal}}
Any cycle appearing in the Chow quotient $Y \chow T$ is a connected,
nodal curve of arithmetic genus 0 and degree $u$, the union of orbit
closures of $U_1$, \dots, $U_i$ where
\beqas
U_j & \neq & (U_j \cap V) \oplus (U_j \cap W), \\
U_1 \cap W & = & 0, \\
U_i \cap V & = & 0, \\
U_j \cap V & = & (U_{j+1} + W) \mod W, \\
\mbox{and }U_{j+1} \cap W & = & (U_j + V) \mod V 
\eeqas 
for all positive $j \leq i-1$. \fp
\es

\hspace{3.6in}
\vspace*{1in}
\font\thinlinefont=cmr5
\begingroup\makeatletter\ifx\SetFigFont\undefined%
\gdef\SetFigFont#1#2#3#4#5{%
  \reset@font\fontsize{#1}{#2pt}%
  \fontfamily{#3}\fontseries{#4}\fontshape{#5}%
  \selectfont}%
\fi\endgroup%
\raisebox{2.5in}{\beginpicture
\setcoordinatesystem units <1.00000cm,1.00000cm>
\unitlength=1.00000cm
\linethickness=1pt
\setplotsymbol ({\makebox(0,0)[l]{\tencirc\symbol{'160}}})
\setshadesymbol ({\thinlinefont .})
\setlinear
%
%
\linethickness= 0.500pt
\setplotsymbol ({\thinlinefont .})
\setdashes < 0.1270cm>
\plot 12.700 24.130 13.335 24.130 /
%
\linethickness= 0.500pt
\setplotsymbol ({\thinlinefont .})
\plot 12.065 23.178 13.970 23.178 /
%
\linethickness= 0.500pt
\setplotsymbol ({\thinlinefont .})
\plot 11.430 22.225 14.605 22.225 /
%
\linethickness= 0.500pt
\setplotsymbol ({\thinlinefont .})
\plot 10.795 21.273 15.240 21.273 /
%
\linethickness= 0.500pt
\setplotsymbol ({\thinlinefont .})
\plot 11.430 20.320 14.605 20.320 /
%
\linethickness= 0.500pt
\setplotsymbol ({\thinlinefont .})
\plot 12.065 19.367 13.970 19.367 /
%
\linethickness= 0.500pt
\setplotsymbol ({\thinlinefont .})
\plot 12.700 18.415 13.335 18.415 /
%
\linethickness= 0.500pt
\setplotsymbol ({\thinlinefont .})
\setsolid
\plot 12.859 24.289 13.335 22.860 /
%
\linethickness= 0.500pt
\setplotsymbol ({\thinlinefont .})
\plot 13.335 23.495 12.859 21.907 /
%
\linethickness= 0.500pt
\setplotsymbol ({\thinlinefont .})
\plot 12.859 22.384 13.494 20.955 /
%
\linethickness= 0.500pt
\setplotsymbol ({\thinlinefont .})
\plot 13.494 21.590 12.859 20.003 /
%
\linethickness= 0.500pt
\setplotsymbol ({\thinlinefont .})
\plot 12.859 20.479 13.494 19.050 /
%
\linethickness= 0.500pt
\setplotsymbol ({\thinlinefont .})
\plot 13.494 19.685 12.859 18.256 /
%
\linethickness= 0.500pt
\setplotsymbol ({\thinlinefont .})
\putrule from 14.605 22.225 to 14.607 22.225
\plot 14.607 22.225 14.616 22.221 /
\plot 14.616 22.221 14.635 22.214 /
\plot 14.635 22.214 14.666 22.204 /
\plot 14.666 22.204 14.711 22.187 /
\plot 14.711 22.187 14.760 22.170 /
\plot 14.760 22.170 14.808 22.151 /
\plot 14.808 22.151 14.855 22.134 /
\plot 14.855 22.134 14.897 22.115 /
\plot 14.897 22.115 14.935 22.100 /
\plot 14.935 22.100 14.969 22.085 /
\plot 14.969 22.085 15.001 22.070 /
\plot 15.001 22.070 15.028 22.056 /
\plot 15.028 22.056 15.056 22.039 /
\plot 15.056 22.039 15.081 22.024 /
\plot 15.081 22.024 15.107 22.007 /
\plot 15.107 22.007 15.132 21.988 /
\plot 15.132 21.988 15.157 21.969 /
\plot 15.157 21.969 15.183 21.948 /
\plot 15.183 21.948 15.210 21.927 /
\plot 15.210 21.927 15.236 21.903 /
\plot 15.236 21.903 15.261 21.880 /
\plot 15.261 21.880 15.287 21.859 /
\plot 15.287 21.859 15.312 21.836 /
\plot 15.312 21.836 15.335 21.812 /
\plot 15.335 21.812 15.356 21.791 /
\plot 15.356 21.791 15.378 21.770 /
\plot 15.378 21.770 15.399 21.749 /
\plot 15.399 21.749 15.420 21.728 /
\plot 15.420 21.728 15.441 21.706 /
\plot 15.441 21.706 15.462 21.685 /
\plot 15.462 21.685 15.483 21.662 /
\plot 15.483 21.662 15.507 21.637 /
\plot 15.507 21.637 15.528 21.611 /
\plot 15.528 21.611 15.551 21.586 /
\plot 15.551 21.586 15.570 21.560 /
\plot 15.570 21.560 15.591 21.533 /
\plot 15.591 21.533 15.608 21.507 /
\plot 15.608 21.507 15.625 21.482 /
\plot 15.625 21.482 15.640 21.457 /
\plot 15.640 21.457 15.653 21.431 /
\plot 15.653 21.431 15.663 21.404 /
\plot 15.663 21.404 15.674 21.378 /
\plot 15.674 21.378 15.682 21.351 /
\plot 15.682 21.351 15.689 21.321 /
\plot 15.689 21.321 15.695 21.292 /
\plot 15.695 21.292 15.699 21.260 /
\plot 15.699 21.260 15.704 21.226 /
\plot 15.704 21.226 15.706 21.192 /
\putrule from 15.706 21.192 to 15.706 21.160
\putrule from 15.706 21.160 to 15.706 21.126
\plot 15.706 21.126 15.704 21.095 /
\plot 15.704 21.095 15.701 21.065 /
\plot 15.701 21.065 15.699 21.035 /
\plot 15.699 21.035 15.695 21.008 /
\plot 15.695 21.008 15.691 20.980 /
\plot 15.691 20.980 15.682 20.951 /
\plot 15.682 20.951 15.676 20.921 /
\plot 15.676 20.921 15.665 20.892 /
\plot 15.665 20.892 15.655 20.862 /
\plot 15.655 20.862 15.642 20.832 /
\plot 15.642 20.832 15.629 20.803 /
\plot 15.629 20.803 15.615 20.775 /
\plot 15.615 20.775 15.600 20.750 /
\plot 15.600 20.750 15.583 20.724 /
\plot 15.583 20.724 15.566 20.703 /
\plot 15.566 20.703 15.549 20.682 /
\plot 15.549 20.682 15.532 20.663 /
\plot 15.532 20.663 15.513 20.646 /
\plot 15.513 20.646 15.492 20.629 /
\plot 15.492 20.629 15.471 20.612 /
\plot 15.471 20.612 15.445 20.595 /
\plot 15.445 20.595 15.420 20.578 /
\plot 15.420 20.578 15.392 20.561 /
\plot 15.392 20.561 15.363 20.546 /
\plot 15.363 20.546 15.333 20.532 /
\plot 15.333 20.532 15.303 20.517 /
\plot 15.303 20.517 15.274 20.504 /
\plot 15.274 20.504 15.244 20.491 /
\plot 15.244 20.491 15.215 20.479 /
\plot 15.215 20.479 15.187 20.468 /
\plot 15.187 20.468 15.160 20.458 /
\plot 15.160 20.458 15.130 20.447 /
\plot 15.130 20.447 15.100 20.436 /
\plot 15.100 20.436 15.069 20.426 /
\plot 15.069 20.426 15.037 20.413 /
\plot 15.037 20.413 15.005 20.403 /
\plot 15.005 20.403 14.973 20.392 /
\plot 14.973 20.392 14.944 20.384 /
\plot 14.944 20.384 14.914 20.373 /
\plot 14.914 20.373 14.887 20.367 /
\plot 14.887 20.367 14.861 20.358 /
\plot 14.861 20.358 14.838 20.352 /
\plot 14.838 20.352 14.817 20.345 /
\plot 14.817 20.345 14.789 20.339 /
\plot 14.789 20.339 14.764 20.335 /
\plot 14.764 20.335 14.738 20.331 /
\plot 14.738 20.331 14.711 20.326 /
\plot 14.711 20.326 14.683 20.324 /
\plot 14.683 20.324 14.656 20.322 /
\putrule from 14.656 20.322 to 14.630 20.322
\plot 14.630 20.322 14.613 20.320 /
\putrule from 14.613 20.320 to 14.607 20.320
\putrule from 14.607 20.320 to 14.605 20.320
%
\linethickness= 0.500pt
\setplotsymbol ({\thinlinefont .})
\putrule from 12.541 24.130 to 12.535 24.130
\putrule from 12.535 24.130 to 12.520 24.130
\putrule from 12.520 24.130 to 12.495 24.130
\plot 12.495 24.130 12.463 24.128 /
\putrule from 12.463 24.128 to 12.425 24.128
\plot 12.425 24.128 12.385 24.126 /
\plot 12.385 24.126 12.344 24.124 /
\plot 12.344 24.124 12.306 24.122 /
\plot 12.306 24.122 12.268 24.119 /
\plot 12.268 24.119 12.232 24.115 /
\plot 12.232 24.115 12.196 24.113 /
\plot 12.196 24.113 12.158 24.109 /
\plot 12.158 24.109 12.118 24.102 /
\plot 12.118 24.102 12.084 24.098 /
\plot 12.084 24.098 12.050 24.094 /
\plot 12.050 24.094 12.012 24.090 /
\plot 12.012 24.090 11.972 24.083 /
\plot 11.972 24.083 11.927 24.075 /
\plot 11.927 24.075 11.883 24.069 /
\plot 11.883 24.069 11.836 24.060 /
\plot 11.836 24.060 11.788 24.050 /
\plot 11.788 24.050 11.737 24.039 /
\plot 11.737 24.039 11.686 24.028 /
\plot 11.686 24.028 11.633 24.016 /
\plot 11.633 24.016 11.582 24.005 /
\plot 11.582 24.005 11.532 23.990 /
\plot 11.532 23.990 11.483 23.978 /
\plot 11.483 23.978 11.434 23.963 /
\plot 11.434 23.963 11.388 23.950 /
\plot 11.388 23.950 11.343 23.933 /
\plot 11.343 23.933 11.299 23.918 /
\plot 11.299 23.918 11.259 23.904 /
\plot 11.259 23.904 11.218 23.887 /
\plot 11.218 23.887 11.178 23.870 /
\plot 11.178 23.870 11.136 23.853 /
\plot 11.136 23.853 11.096 23.834 /
\plot 11.096 23.834 11.053 23.812 /
\plot 11.053 23.812 11.009 23.789 /
\plot 11.009 23.789 10.966 23.766 /
\plot 10.966 23.766 10.922 23.741 /
\plot 10.922 23.741 10.878 23.715 /
\plot 10.878 23.715 10.835 23.688 /
\plot 10.835 23.688 10.791 23.660 /
\plot 10.791 23.660 10.748 23.630 /
\plot 10.748 23.630 10.706 23.601 /
\plot 10.706 23.601 10.666 23.571 /
\plot 10.666 23.571 10.628 23.542 /
\plot 10.628 23.542 10.588 23.510 /
\plot 10.588 23.510 10.552 23.480 /
\plot 10.552 23.480 10.513 23.448 /
\plot 10.513 23.448 10.478 23.415 /
\plot 10.478 23.415 10.448 23.387 /
\plot 10.448 23.387 10.416 23.360 /
\plot 10.416 23.360 10.386 23.330 /
\plot 10.386 23.330 10.355 23.298 /
\plot 10.355 23.298 10.323 23.264 /
\plot 10.323 23.264 10.291 23.228 /
\plot 10.291 23.228 10.259 23.190 /
\plot 10.259 23.190 10.228 23.152 /
\plot 10.228 23.152 10.194 23.110 /
\plot 10.194 23.110 10.162 23.065 /
\plot 10.162 23.065 10.130 23.019 /
\plot 10.130 23.019 10.099 22.970 /
\plot 10.099 22.970 10.069 22.921 /
\plot 10.069 22.921 10.039 22.868 /
\plot 10.039 22.868 10.010 22.816 /
\plot 10.010 22.816  9.982 22.761 /
\plot  9.982 22.761  9.957 22.705 /
\plot  9.957 22.705  9.931 22.648 /
\plot  9.931 22.648  9.908 22.589 /
\plot  9.908 22.589  9.887 22.530 /
\plot  9.887 22.530  9.868 22.468 /
\plot  9.868 22.468  9.849 22.407 /
\plot  9.849 22.407  9.832 22.344 /
\plot  9.832 22.344  9.817 22.278 /
\plot  9.817 22.278  9.804 22.223 /
\plot  9.804 22.223  9.794 22.168 /
\plot  9.794 22.168  9.783 22.111 /
\plot  9.783 22.111  9.773 22.049 /
\plot  9.773 22.049  9.764 21.988 /
\plot  9.764 21.988  9.756 21.924 /
\plot  9.756 21.924  9.749 21.859 /
\plot  9.749 21.859  9.743 21.789 /
\plot  9.743 21.789  9.737 21.719 /
\plot  9.737 21.719  9.732 21.647 /
\plot  9.732 21.647  9.728 21.575 /
\plot  9.728 21.575  9.724 21.501 /
\plot  9.724 21.501  9.722 21.425 /
\putrule from  9.722 21.425 to  9.722 21.349
\plot  9.722 21.349  9.720 21.273 /
\plot  9.720 21.273  9.722 21.196 /
\putrule from  9.722 21.196 to  9.722 21.120
\plot  9.722 21.120  9.724 21.044 /
\plot  9.724 21.044  9.728 20.970 /
\plot  9.728 20.970  9.732 20.898 /
\plot  9.732 20.898  9.737 20.826 /
\plot  9.737 20.826  9.743 20.756 /
\plot  9.743 20.756  9.749 20.686 /
\plot  9.749 20.686  9.756 20.621 /
\plot  9.756 20.621  9.764 20.557 /
\plot  9.764 20.557  9.773 20.496 /
\plot  9.773 20.496  9.783 20.434 /
\plot  9.783 20.434  9.794 20.377 /
\plot  9.794 20.377  9.804 20.322 /
\plot  9.804 20.322  9.817 20.267 /
\plot  9.817 20.267  9.832 20.201 /
\plot  9.832 20.201  9.849 20.138 /
\plot  9.849 20.138  9.868 20.077 /
\plot  9.868 20.077  9.887 20.015 /
\plot  9.887 20.015  9.908 19.956 /
\plot  9.908 19.956  9.931 19.897 /
\plot  9.931 19.897  9.957 19.840 /
\plot  9.957 19.840  9.982 19.784 /
\plot  9.982 19.784 10.010 19.729 /
\plot 10.010 19.729 10.039 19.677 /
\plot 10.039 19.677 10.069 19.624 /
\plot 10.069 19.624 10.099 19.575 /
\plot 10.099 19.575 10.130 19.526 /
\plot 10.130 19.526 10.162 19.480 /
\plot 10.162 19.480 10.194 19.435 /
\plot 10.194 19.435 10.228 19.393 /
\plot 10.228 19.393 10.259 19.355 /
\plot 10.259 19.355 10.291 19.317 /
\plot 10.291 19.317 10.323 19.281 /
\plot 10.323 19.281 10.355 19.247 /
\plot 10.355 19.247 10.386 19.215 /
\plot 10.386 19.215 10.416 19.185 /
\plot 10.416 19.185 10.448 19.158 /
\plot 10.448 19.158 10.478 19.128 /
\plot 10.478 19.128 10.513 19.097 /
\plot 10.513 19.097 10.552 19.065 /
\plot 10.552 19.065 10.588 19.035 /
\plot 10.588 19.035 10.628 19.003 /
\plot 10.628 19.003 10.666 18.974 /
\plot 10.666 18.974 10.706 18.944 /
\plot 10.706 18.944 10.748 18.915 /
\plot 10.748 18.915 10.791 18.885 /
\plot 10.791 18.885 10.835 18.857 /
\plot 10.835 18.857 10.878 18.830 /
\plot 10.878 18.830 10.922 18.804 /
\plot 10.922 18.804 10.966 18.779 /
\plot 10.966 18.779 11.009 18.756 /
\plot 11.009 18.756 11.053 18.733 /
\plot 11.053 18.733 11.096 18.711 /
\plot 11.096 18.711 11.136 18.692 /
\plot 11.136 18.692 11.178 18.675 /
\plot 11.178 18.675 11.218 18.658 /
\plot 11.218 18.658 11.259 18.641 /
\plot 11.259 18.641 11.299 18.627 /
\plot 11.299 18.627 11.343 18.612 /
\plot 11.343 18.612 11.388 18.595 /
\plot 11.388 18.595 11.434 18.582 /
\plot 11.434 18.582 11.483 18.567 /
\plot 11.483 18.567 11.532 18.555 /
\plot 11.532 18.555 11.582 18.540 /
\plot 11.582 18.540 11.633 18.529 /
\plot 11.633 18.529 11.686 18.517 /
\plot 11.686 18.517 11.737 18.506 /
\plot 11.737 18.506 11.788 18.495 /
\plot 11.788 18.495 11.836 18.485 /
\plot 11.836 18.485 11.883 18.476 /
\plot 11.883 18.476 11.927 18.470 /
\plot 11.927 18.470 11.972 18.462 /
\plot 11.972 18.462 12.012 18.455 /
\plot 12.012 18.455 12.050 18.451 /
\plot 12.050 18.451 12.084 18.447 /
\plot 12.084 18.447 12.118 18.440 /
\plot 12.118 18.440 12.158 18.436 /
\plot 12.158 18.436 12.196 18.432 /
\plot 12.196 18.432 12.232 18.430 /
\plot 12.232 18.430 12.268 18.426 /
\plot 12.268 18.426 12.306 18.423 /
\plot 12.306 18.423 12.344 18.421 /
\plot 12.344 18.421 12.385 18.419 /
\plot 12.385 18.419 12.425 18.417 /
\putrule from 12.425 18.417 to 12.463 18.417
\plot 12.463 18.417 12.495 18.415 /
\putrule from 12.495 18.415 to 12.520 18.415
\putrule from 12.520 18.415 to 12.535 18.415
\putrule from 12.535 18.415 to 12.541 18.415
%
\linethickness= 0.500pt
\setplotsymbol ({\thinlinefont .})
\plot 14.764 20.320 14.770 20.314 /
\plot 14.770 20.314 14.783 20.299 /
\plot 14.783 20.299 14.804 20.278 /
\plot 14.804 20.278 14.831 20.246 /
\plot 14.831 20.246 14.861 20.212 /
\plot 14.861 20.212 14.893 20.176 /
\plot 14.893 20.176 14.922 20.142 /
\plot 14.922 20.142 14.948 20.110 /
\plot 14.948 20.110 14.969 20.081 /
\plot 14.969 20.081 14.988 20.053 /
\plot 14.988 20.053 15.005 20.028 /
\plot 15.005 20.028 15.018 20.003 /
\plot 15.018 20.003 15.028 19.975 /
\plot 15.028 19.975 15.039 19.950 /
\plot 15.039 19.950 15.047 19.922 /
\plot 15.047 19.922 15.054 19.892 /
\plot 15.054 19.892 15.060 19.863 /
\plot 15.060 19.863 15.064 19.829 /
\plot 15.064 19.829 15.069 19.795 /
\plot 15.069 19.795 15.071 19.761 /
\putrule from 15.071 19.761 to 15.071 19.725
\putrule from 15.071 19.725 to 15.071 19.691
\plot 15.071 19.691 15.069 19.655 /
\plot 15.069 19.655 15.066 19.622 /
\plot 15.066 19.622 15.064 19.590 /
\plot 15.064 19.590 15.060 19.558 /
\plot 15.060 19.558 15.056 19.526 /
\plot 15.056 19.526 15.050 19.494 /
\plot 15.050 19.494 15.043 19.463 /
\plot 15.043 19.463 15.035 19.431 /
\plot 15.035 19.431 15.026 19.395 /
\plot 15.026 19.395 15.014 19.361 /
\plot 15.014 19.361 15.001 19.325 /
\plot 15.001 19.325 14.988 19.287 /
\plot 14.988 19.287 14.971 19.251 /
\plot 14.971 19.251 14.954 19.215 /
\plot 14.954 19.215 14.933 19.181 /
\plot 14.933 19.181 14.914 19.145 /
\plot 14.914 19.145 14.891 19.113 /
\plot 14.891 19.113 14.867 19.082 /
\plot 14.867 19.082 14.844 19.050 /
\plot 14.844 19.050 14.823 19.027 /
\plot 14.823 19.027 14.800 19.001 /
\plot 14.800 19.001 14.776 18.976 /
\plot 14.776 18.976 14.749 18.951 /
\plot 14.749 18.951 14.721 18.925 /
\plot 14.721 18.925 14.692 18.898 /
\plot 14.692 18.898 14.658 18.872 /
\plot 14.658 18.872 14.626 18.845 /
\plot 14.626 18.845 14.590 18.819 /
\plot 14.590 18.819 14.554 18.794 /
\plot 14.554 18.794 14.518 18.768 /
\plot 14.518 18.768 14.480 18.745 /
\plot 14.480 18.745 14.444 18.722 /
\plot 14.444 18.722 14.406 18.701 /
\plot 14.406 18.701 14.370 18.680 /
\plot 14.370 18.680 14.334 18.661 /
\plot 14.334 18.661 14.298 18.644 /
\plot 14.298 18.644 14.262 18.627 /
\plot 14.262 18.627 14.224 18.612 /
\plot 14.224 18.612 14.186 18.595 /
\plot 14.186 18.595 14.148 18.582 /
\plot 14.148 18.582 14.105 18.567 /
\plot 14.105 18.567 14.059 18.553 /
\plot 14.059 18.553 14.010 18.540 /
\plot 14.010 18.540 13.955 18.525 /
\plot 13.955 18.525 13.898 18.508 /
\plot 13.898 18.508 13.835 18.493 /
\plot 13.835 18.493 13.771 18.479 /
\plot 13.771 18.479 13.705 18.462 /
\plot 13.705 18.462 13.644 18.449 /
\plot 13.644 18.449 13.591 18.436 /
\plot 13.591 18.436 13.549 18.428 /
\plot 13.549 18.428 13.519 18.421 /
\plot 13.519 18.421 13.502 18.417 /
\plot 13.502 18.417 13.496 18.415 /
\putrule from 13.496 18.415 to 13.494 18.415
%
\linethickness= 0.500pt
\setplotsymbol ({\thinlinefont .})
\putrule from 13.494 24.130 to 13.496 24.130
\putrule from 13.496 24.130 to 13.502 24.130
\putrule from 13.502 24.130 to 13.521 24.130
\plot 13.521 24.130 13.553 24.128 /
\putrule from 13.553 24.128 to 13.597 24.128
\plot 13.597 24.128 13.650 24.126 /
\plot 13.650 24.126 13.708 24.124 /
\plot 13.708 24.124 13.767 24.122 /
\plot 13.767 24.122 13.822 24.117 /
\plot 13.822 24.117 13.873 24.113 /
\plot 13.873 24.113 13.919 24.109 /
\plot 13.919 24.109 13.964 24.105 /
\plot 13.964 24.105 14.002 24.100 /
\plot 14.002 24.100 14.036 24.094 /
\plot 14.036 24.094 14.069 24.086 /
\plot 14.069 24.086 14.103 24.077 /
\plot 14.103 24.077 14.133 24.067 /
\plot 14.133 24.067 14.165 24.056 /
\plot 14.165 24.056 14.196 24.043 /
\plot 14.196 24.043 14.230 24.031 /
\plot 14.230 24.031 14.262 24.016 /
\plot 14.262 24.016 14.294 23.999 /
\plot 14.294 23.999 14.328 23.980 /
\plot 14.328 23.980 14.359 23.961 /
\plot 14.359 23.961 14.391 23.939 /
\plot 14.391 23.939 14.421 23.918 /
\plot 14.421 23.918 14.450 23.897 /
\plot 14.450 23.897 14.480 23.876 /
\plot 14.480 23.876 14.506 23.853 /
\plot 14.506 23.853 14.531 23.832 /
\plot 14.531 23.832 14.556 23.808 /
\plot 14.556 23.808 14.580 23.785 /
\plot 14.580 23.785 14.601 23.764 /
\plot 14.601 23.764 14.624 23.738 /
\plot 14.624 23.738 14.645 23.713 /
\plot 14.645 23.713 14.668 23.688 /
\plot 14.668 23.688 14.690 23.658 /
\plot 14.690 23.658 14.711 23.628 /
\plot 14.711 23.628 14.732 23.599 /
\plot 14.732 23.599 14.753 23.567 /
\plot 14.753 23.567 14.772 23.535 /
\plot 14.772 23.535 14.791 23.501 /
\plot 14.791 23.501 14.808 23.470 /
\plot 14.808 23.470 14.823 23.438 /
\plot 14.823 23.438 14.836 23.404 /
\plot 14.836 23.404 14.848 23.372 /
\plot 14.848 23.372 14.859 23.340 /
\plot 14.859 23.340 14.870 23.309 /
\plot 14.870 23.309 14.878 23.277 /
\plot 14.878 23.277 14.887 23.245 /
\plot 14.887 23.245 14.893 23.209 /
\plot 14.893 23.209 14.897 23.175 /
\plot 14.897 23.175 14.903 23.137 /
\plot 14.903 23.137 14.908 23.099 /
\plot 14.908 23.099 14.910 23.061 /
\plot 14.910 23.061 14.912 23.023 /
\putrule from 14.912 23.023 to 14.912 22.985
\putrule from 14.912 22.985 to 14.912 22.947
\putrule from 14.912 22.947 to 14.912 22.911
\plot 14.912 22.911 14.910 22.875 /
\plot 14.910 22.875 14.908 22.843 /
\plot 14.908 22.843 14.903 22.811 /
\plot 14.903 22.811 14.901 22.782 /
\plot 14.901 22.782 14.897 22.754 /
\plot 14.897 22.754 14.889 22.718 /
\plot 14.889 22.718 14.882 22.684 /
\plot 14.882 22.684 14.874 22.650 /
\plot 14.874 22.650 14.863 22.617 /
\plot 14.863 22.617 14.851 22.585 /
\plot 14.851 22.585 14.840 22.553 /
\plot 14.840 22.553 14.827 22.523 /
\plot 14.827 22.523 14.815 22.496 /
\plot 14.815 22.496 14.802 22.473 /
\plot 14.802 22.473 14.789 22.449 /
\plot 14.789 22.449 14.776 22.428 /
\plot 14.776 22.428 14.764 22.409 /
\plot 14.764 22.409 14.749 22.388 /
\plot 14.749 22.388 14.734 22.369 /
\plot 14.734 22.369 14.717 22.348 /
\plot 14.717 22.348 14.696 22.324 /
\plot 14.696 22.324 14.675 22.299 /
\plot 14.675 22.299 14.649 22.274 /
\plot 14.649 22.274 14.628 22.250 /
\plot 14.628 22.250 14.613 22.233 /
\plot 14.613 22.233 14.607 22.227 /
\plot 14.607 22.227 14.605 22.225 /
\linethickness=0pt
\putrectangle corners at  9.694 24.314 and 15.731 18.231
\endpicture}

\vspace*{-6.1in}

\br{Remarks\label{remarks}} 
{\advance \rightskip by200pt 
(1) Of the five equations above, the first
is a nondegeneracy condition, the next two say that the curve
stretches all the way from one extreme fixed component to the other,
and the last two say that the sink of each $T$-orbit is the source of
the next.  Three curves satisfying these conditions are illustrated in
the diagram at right.  The dashed lines represent the components of the
fixed-point set.  The whole picture is reminiscent of the ``broken
flows'' in Morse theory---an analogy which is expanded upon in
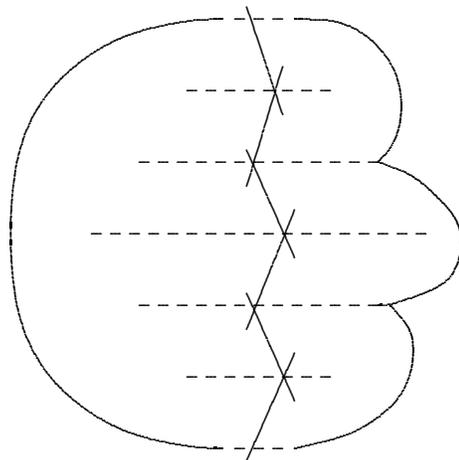
\S\ref{symplectic}.

(2) Each point in the Chow quotient
determines flags in $V$ and $W$ by intersecting those spaces with 

\noindent
the $U_i$.  The {\em Halphen locus} (cf.\ Laksov \cite{bowdoin}) is
the locus where these are maximal.  Equivalently, there are $u$
distinct subspaces $U_1, \dots, U_u$ giving rise to a cycle which is
the union of $u$ lines.  This is illustrated by the middle curve in
the diagram. }

(3) Finally, the alert reader will notice a direct connection with
complete collineations: namely, $U_j$ are the graphs of the linear
maps $f_j$.  Thus the Halphen locus is the locus where every $f_j$ has
rank 1. \er

\br{Proof that \re{gr}{b} = \re{gr}{c}\label{gr-chow/mum}}
Let $\Pj^N = \Pj \La^u (V \oplus W)$.  The Pl\"ucker embedding takes
$Y$ into $\Pj^N$ $T$-equivariantly, and the linearizations $L(\si)$ on
$Y$ extend naturally over $\Pj^N$.  Since the extremal weights of the
$T$-action on $\Pj^N$ are $\pm u$, it follows that the Mumford
quotients $\Pj^N_\si$ are nonempty exactly when $\si \in [0,u]$: the
same range as for $Y$ itself.  Hence the inverse limit of the nonempty
$X_\si$ embeds in that of the $\Pj^N_\si$.

On the other hand, it is easy to check that a generic $T$-orbit
closure in $\Pj^N$ has degree $u$, simply because its source and sink
are in $\Pj \La^u V$ and $\Pj \La^u W$ respectively.  Hence there is a
natural embedding of Chow quotients $Y \chow T \subset \Pj^N \chow T$.

Now the work of Kapranov, Sturmfels, and Zelevinsky \cite{ksz} on
toric varieties shows that the Chow quotient $\Pj^N \chow T$ is an
irreducible component of the inverse limit of Mumford quotients, with
the embedding given by Corollary \ref{corollary}.  Hence the Chow
quotient of $Y$ embeds in the inverse limit of its Mumford quotients.
But both are varieties of the same dimension, and both are
irreducible: the former by definition, the latter by
\reb{vain/gr-mum}. \fp \er

\br{Proof that \re{gr}{b} = \re{gr}{a}\label{gr-chow/hilb}}
By Proposition \ref{nodal} the Chow quotient parametrizes a family of
cycles which are all connected curves of the same arithmetic genus and
degree.  Regarded as reduced schemes, they therefore all have the same
Hilbert polynomial.  On the other hand, the Hilbert quotient
parametrizes a flat family of schemes, and maps to the Chow quotient
by the forgetful morphism taking schemes to cycles.  The generic fiber
of this family is reduced by definition.

But it follows that every fiber is reduced.  Otherwise, it would have
a greater Hilbert polynomial than its reduced subscheme.  Then the
Hilbert polynomial would not be constant in the fibers of the flat
family, which is a contradiction \cite[III 9.9]{har}.

The morphism from the Hilbert quotient to the Chow quotient is
therefore injective as well as surjective.  Since the Chow quotient is
smooth by \reb{vain/gr-mum} and \reb{gr-chow/mum}, Zariski's main theorem
implies that this is an isomorphism.  \fp

\bit{The locus of graphs \label{graphs}}

This section completes the proof of the Main Theorem by relating the
Chow quotient $Y \chow T$ to the locus of graphs of linear maps in
$\Pj V \times \Pj W$.  

Let $\cC$ be the Chow variety parametrizing cycles deformation
equivalent to the graph in $\Pj V \times \Pj W$ of an injective linear
map from a $u$-dimensional subspace of $V$ to $W$.  Then let $\cC'$ be
the closure in $\cC$ of the locus of such graphs.  Similarly, let
$\cH$ be the appropriate Hilbert scheme and $\cH'$ the closure of the
locus of graphs.  The object of this section is to identify both
$\cC'$ and $\cH'$ with the space of complete collineations.

\br{Proof that \re{gr}{b} = \re{graph}{b}\label{gr/graph}}
Consider the diagram
$$\begin{array}{ccc}
\Pj S_{V \oplus W} & \lrow & \Gr_u (V \oplus W) \\
\downarg{} & & \\
\Pj (V \oplus W) & \dasharrow & \Pj V \times \Pj W.
\end{array}$$

The arrow in the bottom row is a rational
map, the geometric quotient of an open set by the $T$-action.  

As in \reb{chow-gr/proj}, the universal cycle on $(Y \chow T) \times
Y$ can be pulled back, pushed forward and divided by the $T$-action to
produce a cycle on $(Y \chow T) \times \Pj V \times \Pj W$.  It is
easy to see that this cycle satisfies Lemma \ref{useful} and hence
defines a morphism from $Y \chow T$ to $\cC$.  The fiber over a point
of $Y \chow T$ defined by subspaces $U_i$ as in Proposition
\ref{nodal} is taken to the cycle
$$
\{ [v],[w] \in \Pj V \times \Pj W \st \la v \oplus \mu w \in U_i
\mbox{ for some } i \mbox{ and some } \la, \mu \in \cx \}.
$$
In particular, a general orbit closure in $Y$ is taken to the graph
of a general linear map $V \to W$.  Furthermore, modulo the action of
$T$, the $U_i$ can be recovered from the above set, so the morphism is
injective.  It therefore remains to show that the derivative is
injective.

As seen in Proposition \ref{nodal}, every cycle appearing in the Chow
quotient $Y \chow T$ is a nodal curve $C$.  Its normal sheaf $N$ is
therefore locally free, and the Zariski tangent space to the Chow
variety is $H^0(C,N)$.  Hence every nonzero deformation in the Chow
quotient is represented by a $T$-invariant normal vector field,
nonzero at the generic point of some component of $C$.

Such a normal vector field pulls back to the inverse image cycle on
$\Pj S_{V\oplus W}$, pushes forward to $\Pj V \times \Pj W$---at least
on the smooth points downstairs---and descends to the quotient by $T$.
On an open subset of the smooth points, it is not killed by the
push-forward, because the tangent vectors to $Y$ killed by projection
from a point $(U, p \in U) \in \Pj S_{V\oplus W}$ are only those
elements of $T_U Y = \Hom\lp U, (V \oplus W)/U\rp$ annihilating $p$.

Thus we obtain a nonzero normal vector field on the smooth points of
the cycle $Z$ in $\Pj V \times \Pj W$ corresponding to $C$.  This
determines an injective map 
$$T_C (Y \chow T) \to H^0(Z_{sm},N_{Z_{sm}/\Pj V \times \Pj W}),$$ 
where $Z_{sm}$ denotes the
smooth part of $Z$.  

On the other hand, there is a natural map $T_Z \cC \to
H^0(Z_{sm},N_{Z_{sm}/\Pj V \times \Pj W})$ which is compatible with
the map described above.  The derivative $T_C (Y \chow T) \to T_Z
\cC$ is therefore a factor of an injective map, so it is injective.
\fp \er

\br{Proof that \re{graph}{b} = \re{graph}{a}\label{graph-chow/hilb}} 
The general element of $\cC'$ is a single orbit closure, so the
natural morphism $\cH' \to \cC'$ is birational and projective.  Now
\reb{gr/graph} has already identified $\cC'$ with the complete
collineations.  An element in the Halphen locus corresponds to
$$\bigcup_{i=1}^u \, \Pj (U_i + W) \mod W  \times \Pj (U_i + V) \mod V \,
 =  \,\bigcup_{i=1}^u \,\Pj^{u-1} \times \Pj^{i-1}.
$$
Regarding this as a reduced scheme, restrict the line bundle
$\co(k,k)$ from $\Pj V \times \Pj W$ to it; the space of sections then
has dimension
$$ \sum_{i=0}^{u-1} \bino{u-2-i+k}{k-1} \bino{i+k}{k}.$$

This is proved by induction on $i$, starting with $\Pj^{u-1}$ and
gluing on each remaining component in turn.  The component $\Pj^{u-i}
\times \Pj^{i-1}$ contributes 
$$\bino{u-2-i+k}{k-1} \bino{i+k}{k} 
= \bino{u-1-i+k}{k} \bino{i+k}{k}
- \bino{u-2-i+k}{k} \bino{i+k}{k}$$ 
dimensions to the space of sections, the first term on the right being
the dimension of its own sections, and the second term imposing the
condition required for sections to agree on the intersection
$\Pj^{u-i} \times \Pj^{i-2}$.

The general fiber, on the other hand, is just a projective space
$\Pj^{u-1}$ projecting linearly to both $\Pj V$ and $\Pj W$.  The
restriction of $\co(k,k)$ is simply $\co(2k)$, so the space of
sections has dimension 
$$\dim H^0\lp\Pj^{u-1}, \co (2k)\rp = \bino{u-1+2k}{2k}.$$  
This turns out to be exactly the same as the summation
above!  One can prove it by multiplying both expressions by the formal
variable $x^{u-1}$ and summing over $u$; after straightforward
manipulations using the identity
$$\sum_{r=0}^\infty \bino{r}{j} x^j = \frac{x^j}{(1-x)^{j+1}},$$
both generating functions can be identified as $(1-x)^{-1-2k}$.  (This
is an application of the ``Snake Oil'' method peddled by Wilf
\cite{wilf}.)

The schemes over the Halphen locus in $\cH'$ must therefore be
reduced.  Otherwise, their Hilbert polynomials would exceed those of
their reduced subschemes, which equal that of the general fiber by the
above.  This would contradict flatness \cite[III 9.9]{har}.  

The fiber in $\cH'$ over each point in the Halphen locus therefore
contains a single closed point.  In particular, it is 0-dimensional.

On the other hand, given any point at all in $Y \chow T$, it is easy
to find a 1-parameter subgroup of $\SL{V} \times \SL{W}$ taking it to
a limit in the Halphen locus.  Just choose any bases for $V$ and $W$
compatible with the flags of \reb{remarks}, and then let $\cx$ act on
the $i$th basis elements of $V$ and $W$ with weights $i$ and $-i$,
respectively.  It is straightforward to check that the limiting cycle
intersects every component of the fixed-point set, and so must be n
the Halphen locus.  (By the way, this is the only appearance in our
story of the natural $\SL{V} \times \SL{W}$-action on complete
collineations.)

Of course, $\SL{V} \times \SL{W}$ acts compatibly on $\cC'$, $\cH'$,
and their universal families, so each fiber of $\cH' \to \cC'$ is
taken by some $\cx$-action to a limit over the Halphen locus.  But in
a projective morphism, the dimension of the fiber is semi-continuous.
Hence all fibers are 0-dimensional, so the morphism is finite
\cite[III Ex.\ 11.2]{har}.  Since it is also birational and $\cC'$ is
smooth, it is an isomorphism by Zariski's main theorem.  \fp \er

\bit{Relation with Gromov-Witten invariants\label{gw}}

There is still one more interpretation of the space of complete
collineations worth mentioning because of its current interest.  It
relates the complete collineations to stable maps in the sense of
Kontsevich \cite{k}.

\bs{Proposition}
Let $\mo{2}(Y,u)$ be the moduli space of 2-pointed stable maps of
genus 0 and degree $u$ to $Y = \Gr_u(V \oplus W)$, and let $\ev:  
\mo{2}(Y,u) \to Y \times Y$ be the evaluation map.  Then the rank $u$
complete collineations are naturally isomorphic to $\ev^{-1}(\Gr_uV
\times \Gr_uW)$.  \es

\pf.  It follows from Proposition \ref{nodal} and \reb{gr-chow/hilb}
that the Hilbert quotient parametrizes a family of stable maps of
genus 0 and degree $u$, all of which are $T$-invariant.  On the other
hand, it is easy to check that these are the only $T$-invariant maps
in $\ev^{-1}(\Gr_uV \times \Gr_uW)$, which is a smooth projective
variety \cite{fp}.  The fixed-point set of the $T$-action therefore
has only one component, so it must be the whole space.
(Alternatively, a dimension count will do.) \fp

For simplicity consider the case where $\dim V = u = \dim
W$, so that $\Gr_uV$ and $\Gr_uW$ are points.  The proposition then
asserts simply that the fiber of $\ev$ over the point $(V,W)$, and
hence by symmetry over a general point, is isomorphic to the complete
collineations.
 
This result can be used to compute some higher-point Gromov-Witten
invariants of Grassmannians which, to the author's knowledge, cannot
at present be evaluated by other means.  One pulls back cohomology
classes from $Y$ to the universal Chow family, then pushes forward to
$Y \chow T$ and integrates.  For example, the class $c_2(S_{V \oplus
W})$ in $H^4(Y, \Q)$ vanishes when it is pulled back and pushed
forward to $Y \chow T$.  This can be shown by computing its value on
lines in the Halphen locus, which are known \cite{bifet} to generate
$H_2$ of the complete collineations.  It follows that for any classes
$\alpha_1, \dots, \alpha_k \in H^*(\Gr_u \C^{2u}, \Q)$,
$$\langle p,p,c_2,\alpha_1, \dots, \alpha_k \rangle_u = 0.$$
Here $\langle \cdots \rangle_u$ denotes the degree $u$ Gromov-Witten
invariant, $p$ is the Poincar\'e dual of a point, and $c_2$ is the
second Chern class of the tautological bundle, as above.

One could pursue this line of reasoning to evaluate more Gromov-Witten
invariants of Grassmannians, and even of Lagrangian and orthogonal
Grassmannians using the results of \S\S\ref{quadrics} and \ref{skew}.
However, there are severe constraints on the degrees and classes
involved, so the total information available is somewhat limited.

\bit{Appendix: complete quadrics\label{quadrics}}

The choice of two vector spaces $V$ and $W$ has governed everything
done so far.  By specializing to the case $W = V^*$, we discover
additional structure.  Specifically, we can ask the objects of study
to be symmetric or anti-symmetric in $V$.

In fact, historically the study of the symmetric objects, namely the
{\em complete quadrics}, preceded that of complete collineations; it
was initiated by Chasles in the 1860's.  A complete quadric is a
finite sequence of quadrics $Q_i$, where $Q_1$ is a hypersurface in
$\Pj V$, $Q_{i+1}$ is a hypersurface in the singular locus of $Q_i$,
and the last $Q_i$ is smooth.  This makes sense, because the singular
locus of a quadric is always a projective subspace.  

As with complete collineations, the notion of a family is
delicate, but it reduces to the corresponding notion for a complete
collineation if each $Q_{i+1}$ is regarded as a self-adjoint map
$\ker Q_i \to \coker Q_i$, up to a scalar.

The moduli space of complete quadrics can be constructed just as in
\reb{defns}, only substituting
$$\Pj \Sym^2 V^* \dasharrow 
\bigtimes_{i=0}^u \Pj \Sym^2 \La^i V^*$$ 
for the rational map given there.

Everything asserted so far in the paper goes through mutatis mutandis
for complete quadrics. 

\bs{Theorem\label{symm}} The following are isomorphic:

\begin{enumerate}

\item The moduli space of complete quadrics in $\Pj V$;
  
\item The variety defined by blowing up $\Pj \Sym^2 V^*$ along the
  image of the Veronese embedding, and then along the proper
  transforms of each of its secant varieties in turn;

\item The closure {\rm (a)} in the Hilbert scheme, or {\rm (b)} in the
Chow variety of $\Pj V \times \Pj V^*$, of the locus of graphs of
invertible self-adjoint linear maps $V \to V^*$;

\item The {\rm (a)} Chow quotient, or {\rm (b)} distinguished
component of the inverse limit of Mumford quotients $X \mod
\SL{U}$, where 
$U$ is a fixed vector space of the same dimension as $V$, and
$X$ is the closure in 
$\Pj \Hom(U, V) \times \Pj \Hom (U,V^*)$ of the 
locus of $(f,g)$ with $f$ invertible and $f^* g$ self-adjoint;

\item The {\rm (a)} Hilbert quotient, or {\rm (b)} Chow quotient, or
  {\rm (c)} inverse limit of Mumford quotients $\LaGr_u (V \oplus V^*)
  \modd \cx$, where $\LaGr$ denotes the Grassmannian of subspaces
  Lagrangian with respect to the natural symplectic form on $V \oplus
  V^*$, and $\cx$ acts with weight $1$ on $V$ and $-1$ on $V^*$.

\end{enumerate}

\es

\noindent {\em Sketch of proof.}  
Each step of the proof consists either of
exactly imitating the corresponding step in the proof of the Main
Theorem, or of exploiting the embedding of an object named here in its
counterpart from the Main Theorem.  The only subtleties are as
follows.

First, one needs to be sure that the Chow quotient and limit of
Mumford quotients for the Lagrangian Grassmannian embed in their
counterparts for the ordinary Grassmannian.  This is accomplished by
checking that the degree of a generic $T$-orbit closure is the same in
each, as in \reb{gr-chow/mum}, and by verifying that the smaller
quotient is nonempty whenever the larger one is.
 
Second, one needs to know that the secant varieties of the Veronese
embedding are exactly the intersection with $\Pj \Sym^2 V^*$ of
the secant varieties of the Segre embedding.  This follows from their
determinantal description: see for example Harris \cite[Ex.\ 9.19]{h}.

Third, in order to be sure that the whole inverse system of Mumford
quotients of the Lagrangian Grassmannian comes from the ordinary
Grassmannian, one must check that $\Pic(\LaGr) = \Z$.  But this
follows from the Borel-Weil theory, since $\LaGr$ is a quotient of a
semisimple group by a maximal parabolic subgroup.

Finally, one would also like to check that $X$ has Picard group $\Z
\oplus \Z$.  This would imply that the inverse limit of its Mumford
quotients is precisely the space of complete quadrics.  However, $X$
has complicated singularities and understanding the Picard group seems
intractable.  Nevertheless, to conclude only that the {\em
distinguished component} of the inverse limit coincides with the
complete quadrics is easy.  After all, the full inverse system of $X
\modd \SL{U}$ has as a subsystem those quotients by linearizations
which are restricted from $\Pj \Hom(U, V) \times \Pj \Hom (U,V^*)$.
The inverse limit therefore maps naturally to the inverse limit of
this subsystem, which as a subspace of the inverse limit of $\Pj
\Hom(U, V) \times \Pj \Hom (U,V^*) \modd \SL{U}$ is known to be the
complete quadrics.  On the other hand, by Theorem \ref{chow->mum}
there is also a morphism to the distinguished component from the Chow
quotient, which again is known to be the complete quadrics.  The
distinguished component is therefore sandwiched between the complete
quadrics by two birational morphisms whose composite is the identity.
\fp

In fact, Vainsencher, Laksov and others studied the more general notion
of a $u$-complete quadric, defined essentially as a complete quadric
in a $u$-dimensional subspace.  This generalization does not fit
together with the Chow quotient formulation as well as the
rank $u$ complete collineations did, so it will not be discussed here.

\bit{Appendix: complete skew forms\label{skew}}

Although they are not as well documented in the 19th- or 20th-century
literature, anti-symmetric objects can be studied just as well as
symmetric ones.  The basic definition is the following: a {\em
complete skew form} is a finite sequence of nonzero 2-forms
$\omega_i$, where $\omega_1$ is a 2-form on $V$, $\omega_{i+1}$ is a
2-form on the null space of $\omega_i$, and the last $\omega_i$ is
either nondegenerate or has 1-dimensional null space.

The moduli space can be constructed as in \reb{defns}, only
substituting
$$\Pj \La^2 V^* \dasharrow 
\bigtimes_{i=0}^{u-2} \Pj \La^2 \La^i V^*$$ 
for the rational map given there.

With this understood, there is also an anti-symmetric counterpart of
the Main Theorem.

\bs{Theorem} The following are isomorphic:

\begin{enumerate}

\item The moduli space of complete skew forms on $V$;
  
\item The variety defined by blowing up $\Pj \La^2 V^*$ along the
  image of the Pl\"ucker embedding of $\Gr_2 (V^*)$, and then along
  the proper transforms of each of its secant varieties in turn;

\item The closure {\rm (a)} in the Hilbert scheme, or {\rm (b)} in the
Chow variety of $\Pj V \times \Pj V^*$, of either
{\rm (i)} the locus of graphs of
invertible skew-adjoint linear maps $f: V \to V^*$, if $u$ is even, 
or {\rm (ii)} the locus of subvarieties 
$$\{ ([v],[w]) \in \Pj V \times \Pj f(V) \st f(v) = \la w 
\mbox{ for some } \la \in \C \}$$
for $f: V \to V^*$ skew-adjoint of rank $u-1$, if $u$ is odd;

\item The {\rm (a)} Chow quotient, or {\rm (b)} distinguished
component of the inverse limit of Mumford quotients $Z \mod
\SL{U}$, where 
$U$ is a fixed vector space of the same dimension as $V$, and
$Z$ is the closure in 
$\Pj \Hom(U, V) \times \Pj \Hom (U,V^*)$ of the 
locus of $(f,g)$ with $f$ invertible and $f^* g$ skew-adjoint;

\item The {\rm (a)} Hilbert quotient, or {\rm (b)} Chow quotient, or
  {\rm (c)} inverse limit of Mumford quotients $\OGr^+_u (V \oplus
  V^*) \modd \cx$, where $\OGr$ denotes the Grassmannian of subspaces
  orthogonal with respect to the natural symmetric 2-form on $V \oplus
  V^*$, $\OGr^+$ denotes the connected component containing $V$, and
  $\cx$ acts with weight $1$ on $V$ and $-1$ on $V^*$.

\end{enumerate}

\es

\noindent {\em Sketch of proof.}  
As before, the proof mainly parallels that of
the Main Theorem.  However, in addition to all of the subtleties
mentioned under Theorem \ref{symm} above, a few new ones appear.

The most salient point is that the orthogonal Grassmannian is
reducible: it has two connected components, but only one is needed.
We choose the component $\OGr^+$ containing $V \subset V \oplus V^*$.
Any element of $\OGr^+$ has even-dimensional intersection with
$V$, so $\OGr^+$ intersects only {\em every other} component of the
fixed-point set.  Blow-ups and blow-downs of the quotient therefore
appear only when $\si$ crosses an even integer.

The blow-up loci are smooth, and as in Proposition
\ref{blow-up/secant}, the birational map to the first quotient takes
them to the secant varieties of the Pl\"ucker embedding.  They are
therefore exactly the proper transforms of these secant varieties.

Once this is taken into account, everything follows the proof for
complete collineations.  As with complete quadrics, it is better at
some points to follow a parallel but independent argument, at others
to use the embedding of each item in its counterpart from the Main
Theorem.

In the case when $u$ is odd, the embeddings differ slightly from the
previous cases.  For instance, the complete skew forms are not
complete collineations in the obvious way, because there are no
invertible skew-adjoint linear maps.  Nevertheless, they embed in the
complete collineations as follows: to each sequence $\omega_1, \dots,
\omega_k$ of skew forms, associate the corresponding skew-adjoint
linear maps $f_1: V \to V^*$, $f_2 : \ker f_1 \to \coker f_1 = \ker
f_1^*$, etc., and then append one final nonzero map $f_{k+1}: \ker f_k
\to \coker f_k$, whose choice is unique up to a scalar.

Similarly, all the other constructions also embed in their counterparts
from the Main Theorem; it is straightforward to find all these
embeddings and show that they are compatible. \fp

The diagram below portrays typical curves appearing in the Chow
quotient of the orthogonal Grassmannian.  The case $u$ even is shown
on the left; the case $u$ odd, on the right.  The dashed lines
represent the fixed components of the $T$-action on the full
Grassmannian.  Two curves are shown in each case: the first contained
in the good component $\OGr^+$ and the second contained in the other
discarded component.  Notice that they only visit every other fixed
component.  In the even case, the good curve stretches all the way
from one extreme fixed component to the other, so it actually comes
from the Chow quotient of the full Grassmannian.  In the odd case,
even the good curve does not make it all the way, so an extra line
(shown in dots) must be added to define an embedding of $\OGr^+ \chow
T$ in $\Gr \chow T$.
 
\vspace{.5cm}
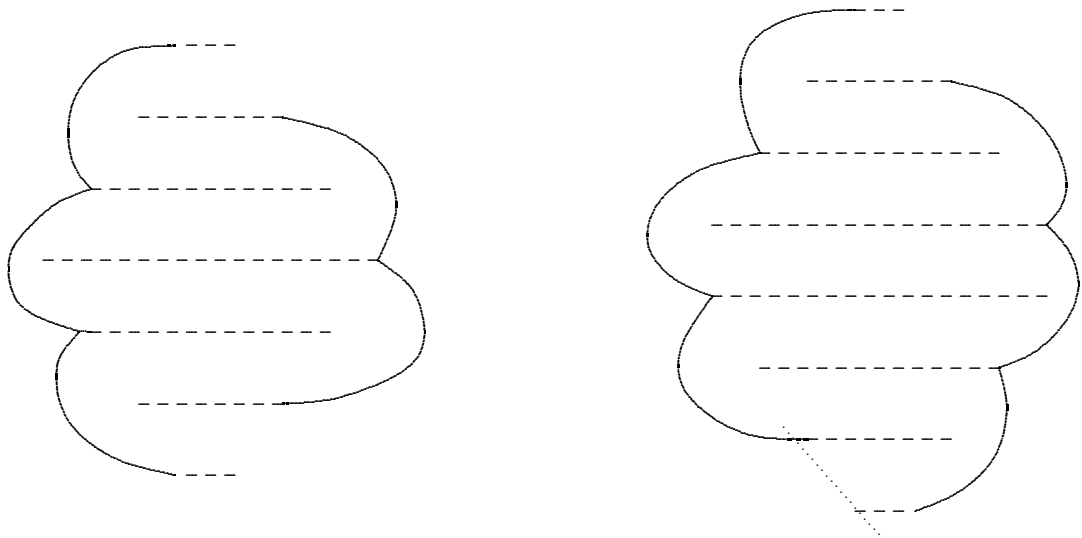
\begin{center}
\font\thinlinefont=cmr5
\begingroup\makeatletter\ifx\SetFigFont\undefined%
\gdef\SetFigFont#1#2#3#4#5{%
  \reset@font\fontsize{#1}{#2pt}%
  \fontfamily{#3}\fontseries{#4}\fontshape{#5}%
  \selectfont}%
\fi\endgroup%
\mbox{\beginpicture
\setcoordinatesystem units <1.00000cm,1.00000cm>
\unitlength=1.00000cm
\linethickness=1pt
\setplotsymbol ({\makebox(0,0)[l]{\tencirc\symbol{'160}}})
\setshadesymbol ({\thinlinefont .})
\setlinear
%
%
\linethickness= 0.500pt
\setplotsymbol ({\thinlinefont .})
\setdashes < 0.1270cm>
\plot 16.510 24.606 15.875 24.606 /
%
\linethickness= 0.500pt
\setplotsymbol ({\thinlinefont .})
\plot 17.145 23.654 15.240 23.654 /
%
\linethickness= 0.500pt
\setplotsymbol ({\thinlinefont .})
\plot 17.780 22.701 14.605 22.701 /
%
\linethickness= 0.500pt
\setplotsymbol ({\thinlinefont .})
\plot 18.415 21.749 13.970 21.749 /
%
\linethickness= 0.500pt
\setplotsymbol ({\thinlinefont .})
\plot 13.970 20.796 18.415 20.796 /
%
\linethickness= 0.500pt
\setplotsymbol ({\thinlinefont .})
\plot 17.780 19.844 14.605 19.844 /
%
\linethickness= 0.500pt
\setplotsymbol ({\thinlinefont .})
\plot 17.145 18.891 15.240 18.891 /
%
\linethickness= 0.500pt
\setplotsymbol ({\thinlinefont .})
\plot 15.875 17.939 16.510 17.939 /
%
\linethickness= 0.500pt
\setplotsymbol ({\thinlinefont .})
\setsolid
\putrule from 15.875 24.606 to 15.873 24.606
\putrule from 15.873 24.606 to 15.867 24.606
\putrule from 15.867 24.606 to 15.847 24.606
\plot 15.847 24.606 15.814 24.604 /
\putrule from 15.814 24.604 to 15.769 24.604
\plot 15.769 24.604 15.714 24.602 /
\plot 15.714 24.602 15.655 24.600 /
\plot 15.655 24.600 15.593 24.598 /
\plot 15.593 24.598 15.534 24.594 /
\plot 15.534 24.594 15.479 24.589 /
\plot 15.479 24.589 15.428 24.585 /
\plot 15.428 24.585 15.382 24.581 /
\plot 15.382 24.581 15.337 24.577 /
\plot 15.337 24.577 15.295 24.570 /
\plot 15.295 24.570 15.255 24.562 /
\plot 15.255 24.562 15.215 24.553 /
\plot 15.215 24.553 15.177 24.545 /
\plot 15.177 24.545 15.138 24.534 /
\plot 15.138 24.534 15.100 24.524 /
\plot 15.100 24.524 15.062 24.513 /
\plot 15.062 24.513 15.020 24.498 /
\plot 15.020 24.498 14.980 24.483 /
\plot 14.980 24.483 14.937 24.469 /
\plot 14.937 24.469 14.895 24.452 /
\plot 14.895 24.452 14.855 24.433 /
\plot 14.855 24.433 14.812 24.414 /
\plot 14.812 24.414 14.774 24.392 /
\plot 14.774 24.392 14.736 24.371 /
\plot 14.736 24.371 14.700 24.350 /
\plot 14.700 24.350 14.666 24.329 /
\plot 14.666 24.329 14.635 24.306 /
\plot 14.635 24.306 14.605 24.282 /
\plot 14.605 24.282 14.577 24.259 /
\plot 14.577 24.259 14.552 24.236 /
\plot 14.552 24.236 14.527 24.208 /
\plot 14.527 24.208 14.501 24.179 /
\plot 14.501 24.179 14.478 24.147 /
\plot 14.478 24.147 14.457 24.113 /
\plot 14.457 24.113 14.438 24.079 /
\plot 14.438 24.079 14.419 24.041 /
\plot 14.419 24.041 14.404 24.003 /
\plot 14.404 24.003 14.389 23.965 /
\plot 14.389 23.965 14.376 23.925 /
\plot 14.376 23.925 14.366 23.882 /
\plot 14.366 23.882 14.357 23.844 /
\plot 14.357 23.844 14.351 23.804 /
\plot 14.351 23.804 14.347 23.766 /
\plot 14.347 23.766 14.343 23.728 /
\plot 14.343 23.728 14.340 23.690 /
\putrule from 14.340 23.690 to 14.340 23.654
\putrule from 14.340 23.654 to 14.340 23.618
\plot 14.340 23.618 14.343 23.580 /
\plot 14.343 23.580 14.345 23.542 /
\plot 14.345 23.542 14.349 23.503 /
\plot 14.349 23.503 14.353 23.463 /
\plot 14.353 23.463 14.359 23.425 /
\plot 14.359 23.425 14.366 23.383 /
\plot 14.366 23.383 14.374 23.343 /
\plot 14.374 23.343 14.383 23.305 /
\plot 14.383 23.305 14.391 23.266 /
\plot 14.391 23.266 14.400 23.228 /
\plot 14.400 23.228 14.410 23.194 /
\plot 14.410 23.194 14.419 23.161 /
\plot 14.419 23.161 14.427 23.129 /
\plot 14.427 23.129 14.438 23.099 /
\plot 14.438 23.099 14.446 23.072 /
\plot 14.446 23.072 14.459 23.036 /
\plot 14.459 23.036 14.472 23.002 /
\plot 14.472 23.002 14.484 22.968 /
\plot 14.484 22.968 14.499 22.930 /
\plot 14.499 22.930 14.516 22.892 /
\plot 14.516 22.892 14.535 22.849 /
\plot 14.535 22.849 14.556 22.805 /
\plot 14.556 22.805 14.575 22.765 /
\plot 14.575 22.765 14.590 22.733 /
\plot 14.590 22.733 14.601 22.712 /
\plot 14.601 22.712 14.605 22.703 /
\putrule from 14.605 22.703 to 14.605 22.701
%
\linethickness= 0.500pt
\setplotsymbol ({\thinlinefont .})
\putrule from 14.605 22.701 to 14.603 22.701
\plot 14.603 22.701 14.597 22.699 /
\plot 14.597 22.699 14.580 22.695 /
\plot 14.580 22.695 14.550 22.689 /
\plot 14.550 22.689 14.508 22.680 /
\plot 14.508 22.680 14.455 22.667 /
\plot 14.455 22.667 14.393 22.655 /
\plot 14.393 22.655 14.328 22.638 /
\plot 14.328 22.638 14.264 22.623 /
\plot 14.264 22.623 14.201 22.608 /
\plot 14.201 22.608 14.144 22.591 /
\plot 14.144 22.591 14.089 22.576 /
\plot 14.089 22.576 14.040 22.564 /
\plot 14.040 22.564 13.993 22.549 /
\plot 13.993 22.549 13.951 22.534 /
\plot 13.951 22.534 13.913 22.521 /
\plot 13.913 22.521 13.875 22.504 /
\plot 13.875 22.504 13.839 22.490 /
\plot 13.839 22.490 13.801 22.473 /
\plot 13.801 22.473 13.765 22.456 /
\plot 13.765 22.456 13.729 22.437 /
\plot 13.729 22.437 13.693 22.416 /
\plot 13.693 22.416 13.655 22.394 /
\plot 13.655 22.394 13.619 22.371 /
\plot 13.619 22.371 13.581 22.346 /
\plot 13.581 22.346 13.545 22.320 /
\plot 13.545 22.320 13.509 22.293 /
\plot 13.509 22.293 13.473 22.265 /
\plot 13.473 22.265 13.441 22.238 /
\plot 13.441 22.238 13.407 22.210 /
\plot 13.407 22.210 13.377 22.181 /
\plot 13.377 22.181 13.350 22.153 /
\plot 13.350 22.153 13.322 22.123 /
\plot 13.322 22.123 13.299 22.096 /
\plot 13.299 22.096 13.276 22.068 /
\plot 13.276 22.068 13.257 22.039 /
\plot 13.257 22.039 13.233 22.007 /
\plot 13.233 22.007 13.212 21.975 /
\plot 13.212 21.975 13.193 21.939 /
\plot 13.193 21.939 13.176 21.903 /
\plot 13.176 21.903 13.159 21.867 /
\plot 13.159 21.867 13.145 21.827 /
\plot 13.145 21.827 13.134 21.789 /
\plot 13.134 21.789 13.123 21.749 /
\plot 13.123 21.749 13.115 21.709 /
\plot 13.115 21.709 13.109 21.670 /
\plot 13.109 21.670 13.104 21.630 /
\putrule from 13.104 21.630 to 13.104 21.594
\plot 13.104 21.594 13.106 21.558 /
\plot 13.106 21.558 13.109 21.522 /
\plot 13.109 21.522 13.115 21.491 /
\plot 13.115 21.491 13.123 21.457 /
\plot 13.123 21.457 13.134 21.427 /
\plot 13.134 21.427 13.147 21.395 /
\plot 13.147 21.395 13.161 21.364 /
\plot 13.161 21.364 13.178 21.332 /
\plot 13.178 21.332 13.197 21.298 /
\plot 13.197 21.298 13.221 21.266 /
\plot 13.221 21.266 13.244 21.234 /
\plot 13.244 21.234 13.269 21.205 /
\plot 13.269 21.205 13.297 21.173 /
\plot 13.297 21.173 13.324 21.146 /
\plot 13.324 21.146 13.354 21.118 /
\plot 13.354 21.118 13.382 21.093 /
\plot 13.382 21.093 13.411 21.069 /
\plot 13.411 21.069 13.439 21.046 /
\plot 13.439 21.046 13.466 21.027 /
\plot 13.466 21.027 13.494 21.008 /
\plot 13.494 21.008 13.526 20.989 /
\plot 13.526 20.989 13.557 20.970 /
\plot 13.557 20.970 13.591 20.951 /
\plot 13.591 20.951 13.627 20.932 /
\plot 13.627 20.932 13.669 20.915 /
\plot 13.669 20.915 13.714 20.894 /
\plot 13.714 20.894 13.763 20.875 /
\plot 13.763 20.875 13.813 20.856 /
\plot 13.813 20.856 13.862 20.836 /
\plot 13.862 20.836 13.906 20.820 /
\plot 13.906 20.820 13.940 20.807 /
\plot 13.940 20.807 13.959 20.800 /
\plot 13.959 20.800 13.968 20.796 /
\putrule from 13.968 20.796 to 13.970 20.796
%
\linethickness= 0.500pt
\setplotsymbol ({\thinlinefont .})
\plot 13.970 20.796 13.966 20.790 /
\plot 13.966 20.790 13.955 20.775 /
\plot 13.955 20.775 13.938 20.752 /
\plot 13.938 20.752 13.913 20.718 /
\plot 13.913 20.718 13.883 20.673 /
\plot 13.883 20.673 13.849 20.625 /
\plot 13.849 20.625 13.811 20.572 /
\plot 13.811 20.572 13.775 20.519 /
\plot 13.775 20.519 13.741 20.468 /
\plot 13.741 20.468 13.710 20.419 /
\plot 13.710 20.419 13.682 20.375 /
\plot 13.682 20.375 13.657 20.335 /
\plot 13.657 20.335 13.636 20.297 /
\plot 13.636 20.297 13.617 20.261 /
\plot 13.617 20.261 13.600 20.227 /
\plot 13.600 20.227 13.587 20.193 /
\plot 13.587 20.193 13.574 20.161 /
\plot 13.574 20.161 13.559 20.125 /
\plot 13.559 20.125 13.549 20.087 /
\plot 13.549 20.087 13.540 20.049 /
\plot 13.540 20.049 13.532 20.011 /
\plot 13.532 20.011 13.526 19.971 /
\plot 13.526 19.971 13.521 19.931 /
\plot 13.521 19.931 13.519 19.888 /
\putrule from 13.519 19.888 to 13.519 19.848
\plot 13.519 19.848 13.521 19.806 /
\plot 13.521 19.806 13.528 19.765 /
\plot 13.528 19.765 13.534 19.727 /
\plot 13.534 19.727 13.542 19.689 /
\plot 13.542 19.689 13.555 19.653 /
\plot 13.555 19.653 13.568 19.617 /
\plot 13.568 19.617 13.583 19.586 /
\plot 13.583 19.586 13.600 19.552 /
\plot 13.600 19.552 13.619 19.522 /
\plot 13.619 19.522 13.640 19.490 /
\plot 13.640 19.490 13.665 19.459 /
\plot 13.665 19.459 13.693 19.427 /
\plot 13.693 19.427 13.722 19.393 /
\plot 13.722 19.393 13.754 19.361 /
\plot 13.754 19.361 13.788 19.329 /
\plot 13.788 19.329 13.822 19.300 /
\plot 13.822 19.300 13.858 19.268 /
\plot 13.858 19.268 13.896 19.241 /
\plot 13.896 19.241 13.932 19.213 /
\plot 13.932 19.213 13.968 19.188 /
\plot 13.968 19.188 14.004 19.164 /
\plot 14.004 19.164 14.038 19.141 /
\plot 14.038 19.141 14.069 19.122 /
\plot 14.069 19.122 14.103 19.103 /
\plot 14.103 19.103 14.133 19.086 /
\plot 14.133 19.086 14.165 19.069 /
\plot 14.165 19.069 14.199 19.052 /
\plot 14.199 19.052 14.230 19.037 /
\plot 14.230 19.037 14.264 19.022 /
\plot 14.264 19.022 14.298 19.008 /
\plot 14.298 19.008 14.332 18.995 /
\plot 14.332 18.995 14.368 18.982 /
\plot 14.368 18.982 14.402 18.970 /
\plot 14.402 18.970 14.436 18.959 /
\plot 14.436 18.959 14.470 18.951 /
\plot 14.470 18.951 14.503 18.942 /
\plot 14.503 18.942 14.535 18.934 /
\plot 14.535 18.934 14.569 18.927 /
\plot 14.569 18.927 14.601 18.923 /
\plot 14.601 18.923 14.633 18.917 /
\plot 14.633 18.917 14.664 18.912 /
\plot 14.664 18.912 14.698 18.910 /
\plot 14.698 18.910 14.732 18.906 /
\plot 14.732 18.906 14.770 18.904 /
\plot 14.770 18.904 14.815 18.902 /
\plot 14.815 18.902 14.861 18.900 /
\plot 14.861 18.900 14.912 18.898 /
\plot 14.912 18.898 14.967 18.895 /
\putrule from 14.967 18.895 to 15.026 18.895
\plot 15.026 18.895 15.083 18.893 /
\putrule from 15.083 18.893 to 15.136 18.893
\plot 15.136 18.893 15.181 18.891 /
\putrule from 15.181 18.891 to 15.212 18.891
\putrule from 15.212 18.891 to 15.232 18.891
\putrule from 15.232 18.891 to 15.238 18.891
\putrule from 15.238 18.891 to 15.240 18.891
%
\linethickness= 0.500pt
\setplotsymbol ({\thinlinefont .})
\putrule from 17.145 23.654 to 17.147 23.654
\plot 17.147 23.654 17.153 23.652 /
\plot 17.153 23.652 17.170 23.647 /
\plot 17.170 23.647 17.200 23.641 /
\plot 17.200 23.641 17.242 23.630 /
\plot 17.242 23.630 17.295 23.618 /
\plot 17.295 23.618 17.357 23.603 /
\plot 17.357 23.603 17.422 23.588 /
\plot 17.422 23.588 17.486 23.571 /
\plot 17.486 23.571 17.549 23.554 /
\plot 17.549 23.554 17.606 23.537 /
\plot 17.606 23.537 17.661 23.520 /
\plot 17.661 23.520 17.710 23.503 /
\plot 17.710 23.503 17.757 23.489 /
\plot 17.757 23.489 17.799 23.472 /
\plot 17.799 23.472 17.837 23.453 /
\plot 17.837 23.453 17.875 23.436 /
\plot 17.875 23.436 17.913 23.415 /
\plot 17.913 23.415 17.949 23.396 /
\plot 17.949 23.396 17.985 23.372 /
\plot 17.985 23.372 18.021 23.349 /
\plot 18.021 23.349 18.057 23.321 /
\plot 18.057 23.321 18.095 23.294 /
\plot 18.095 23.294 18.131 23.264 /
\plot 18.131 23.264 18.169 23.233 /
\plot 18.169 23.233 18.205 23.201 /
\plot 18.205 23.201 18.241 23.167 /
\plot 18.241 23.167 18.277 23.131 /
\plot 18.277 23.131 18.309 23.097 /
\plot 18.309 23.097 18.343 23.061 /
\plot 18.343 23.061 18.373 23.027 /
\plot 18.373 23.027 18.400 22.993 /
\plot 18.400 22.993 18.428 22.957 /
\plot 18.428 22.957 18.451 22.926 /
\plot 18.451 22.926 18.474 22.892 /
\plot 18.474 22.892 18.495 22.860 /
\plot 18.495 22.860 18.517 22.824 /
\plot 18.517 22.824 18.538 22.786 /
\plot 18.538 22.786 18.557 22.748 /
\plot 18.557 22.748 18.574 22.710 /
\plot 18.574 22.710 18.593 22.669 /
\plot 18.593 22.669 18.608 22.631 /
\plot 18.608 22.631 18.622 22.589 /
\plot 18.622 22.589 18.635 22.549 /
\plot 18.635 22.549 18.646 22.511 /
\plot 18.646 22.511 18.656 22.473 /
\plot 18.656 22.473 18.663 22.435 /
\plot 18.663 22.435 18.669 22.401 /
\plot 18.669 22.401 18.673 22.367 /
\plot 18.673 22.367 18.677 22.335 /
\plot 18.677 22.335 18.680 22.305 /
\putrule from 18.680 22.305 to 18.680 22.278
\putrule from 18.680 22.278 to 18.680 22.242
\plot 18.680 22.242 18.675 22.208 /
\plot 18.675 22.208 18.671 22.174 /
\plot 18.671 22.174 18.665 22.140 /
\plot 18.665 22.140 18.656 22.109 /
\plot 18.656 22.109 18.646 22.077 /
\plot 18.646 22.077 18.635 22.047 /
\plot 18.635 22.047 18.622 22.020 /
\plot 18.622 22.020 18.610 21.996 /
\plot 18.610 21.996 18.599 21.973 /
\plot 18.599 21.973 18.586 21.952 /
\plot 18.586 21.952 18.574 21.933 /
\plot 18.574 21.933 18.559 21.912 /
\plot 18.559 21.912 18.544 21.893 /
\plot 18.544 21.893 18.527 21.872 /
\plot 18.527 21.872 18.506 21.848 /
\plot 18.506 21.848 18.485 21.823 /
\plot 18.485 21.823 18.459 21.797 /
\plot 18.459 21.797 18.438 21.774 /
\plot 18.438 21.774 18.423 21.757 /
\plot 18.423 21.757 18.417 21.751 /
\plot 18.417 21.751 18.415 21.749 /
%
\linethickness= 0.500pt
\setplotsymbol ({\thinlinefont .})
\plot 18.415 21.749 18.421 21.742 /
\plot 18.421 21.742 18.434 21.728 /
\plot 18.434 21.728 18.455 21.706 /
\plot 18.455 21.706 18.485 21.675 /
\plot 18.485 21.675 18.517 21.641 /
\plot 18.517 21.641 18.548 21.605 /
\plot 18.548 21.605 18.580 21.571 /
\plot 18.580 21.571 18.608 21.539 /
\plot 18.608 21.539 18.633 21.510 /
\plot 18.633 21.510 18.654 21.482 /
\plot 18.654 21.482 18.673 21.457 /
\plot 18.673 21.457 18.690 21.431 /
\plot 18.690 21.431 18.707 21.404 /
\plot 18.707 21.404 18.722 21.378 /
\plot 18.722 21.378 18.737 21.351 /
\plot 18.737 21.351 18.749 21.321 /
\plot 18.749 21.321 18.764 21.292 /
\plot 18.764 21.292 18.777 21.258 /
\plot 18.777 21.258 18.790 21.224 /
\plot 18.790 21.224 18.800 21.190 /
\plot 18.800 21.190 18.811 21.154 /
\plot 18.811 21.154 18.819 21.120 /
\plot 18.819 21.120 18.826 21.084 /
\plot 18.826 21.084 18.832 21.050 /
\plot 18.832 21.050 18.836 21.018 /
\plot 18.836 21.018 18.838 20.987 /
\putrule from 18.838 20.987 to 18.838 20.955
\putrule from 18.838 20.955 to 18.838 20.923
\plot 18.838 20.923 18.836 20.892 /
\plot 18.836 20.892 18.832 20.860 /
\plot 18.832 20.860 18.826 20.824 /
\plot 18.826 20.824 18.817 20.790 /
\plot 18.817 20.790 18.809 20.754 /
\plot 18.809 20.754 18.798 20.716 /
\plot 18.798 20.716 18.783 20.680 /
\plot 18.783 20.680 18.771 20.644 /
\plot 18.771 20.644 18.754 20.610 /
\plot 18.754 20.610 18.737 20.574 /
\plot 18.737 20.574 18.718 20.542 /
\plot 18.718 20.542 18.701 20.511 /
\plot 18.701 20.511 18.680 20.479 /
\plot 18.680 20.479 18.661 20.451 /
\plot 18.661 20.451 18.641 20.424 /
\plot 18.641 20.424 18.618 20.396 /
\plot 18.618 20.396 18.595 20.367 /
\plot 18.595 20.367 18.570 20.337 /
\plot 18.570 20.337 18.544 20.307 /
\plot 18.544 20.307 18.517 20.278 /
\plot 18.517 20.278 18.487 20.248 /
\plot 18.487 20.248 18.457 20.218 /
\plot 18.457 20.218 18.428 20.191 /
\plot 18.428 20.191 18.398 20.165 /
\plot 18.398 20.165 18.368 20.140 /
\plot 18.368 20.140 18.339 20.117 /
\plot 18.339 20.117 18.311 20.094 /
\plot 18.311 20.094 18.284 20.074 /
\plot 18.284 20.074 18.256 20.055 /
\plot 18.256 20.055 18.224 20.036 /
\plot 18.224 20.036 18.193 20.017 /
\plot 18.193 20.017 18.159 19.998 /
\plot 18.159 19.998 18.123 19.979 /
\plot 18.123 19.979 18.081 19.962 /
\plot 18.081 19.962 18.036 19.941 /
\plot 18.036 19.941 17.987 19.922 /
\plot 17.987 19.922 17.937 19.903 /
\plot 17.937 19.903 17.888 19.884 /
\plot 17.888 19.884 17.843 19.867 /
\plot 17.843 19.867 17.810 19.854 /
\plot 17.810 19.854 17.791 19.848 /
\plot 17.791 19.848 17.782 19.844 /
\putrule from 17.782 19.844 to 17.780 19.844
%
\linethickness= 0.500pt
\setplotsymbol ({\thinlinefont .})
\putrule from 17.780 19.844 to 17.780 19.842
\plot 17.780 19.842 17.782 19.833 /
\plot 17.782 19.833 17.788 19.812 /
\plot 17.788 19.812 17.799 19.778 /
\plot 17.799 19.778 17.812 19.734 /
\plot 17.812 19.734 17.827 19.681 /
\plot 17.827 19.681 17.839 19.628 /
\plot 17.839 19.628 17.852 19.575 /
\plot 17.852 19.575 17.863 19.524 /
\plot 17.863 19.524 17.871 19.480 /
\plot 17.871 19.480 17.877 19.437 /
\plot 17.877 19.437 17.882 19.395 /
\plot 17.882 19.395 17.886 19.355 /
\putrule from 17.886 19.355 to 17.886 19.315
\putrule from 17.886 19.315 to 17.886 19.281
\plot 17.886 19.281 17.884 19.247 /
\plot 17.884 19.247 17.882 19.213 /
\plot 17.882 19.213 17.877 19.175 /
\plot 17.877 19.175 17.873 19.137 /
\plot 17.873 19.137 17.867 19.097 /
\plot 17.867 19.097 17.860 19.054 /
\plot 17.860 19.054 17.852 19.014 /
\plot 17.852 19.014 17.841 18.970 /
\plot 17.841 18.970 17.829 18.927 /
\plot 17.829 18.927 17.818 18.887 /
\plot 17.818 18.887 17.803 18.845 /
\plot 17.803 18.845 17.788 18.804 /
\plot 17.788 18.804 17.774 18.766 /
\plot 17.774 18.766 17.757 18.728 /
\plot 17.757 18.728 17.740 18.694 /
\plot 17.740 18.694 17.721 18.661 /
\plot 17.721 18.661 17.702 18.627 /
\plot 17.702 18.627 17.681 18.595 /
\plot 17.681 18.595 17.657 18.563 /
\plot 17.657 18.563 17.634 18.531 /
\plot 17.634 18.531 17.606 18.502 /
\plot 17.606 18.502 17.579 18.470 /
\plot 17.579 18.470 17.551 18.438 /
\plot 17.551 18.438 17.520 18.409 /
\plot 17.520 18.409 17.488 18.379 /
\plot 17.488 18.379 17.456 18.349 /
\plot 17.456 18.349 17.422 18.322 /
\plot 17.422 18.322 17.388 18.296 /
\plot 17.388 18.296 17.355 18.271 /
\plot 17.355 18.271 17.323 18.246 /
\plot 17.323 18.246 17.291 18.224 /
\plot 17.291 18.224 17.259 18.203 /
\plot 17.259 18.203 17.230 18.184 /
\plot 17.230 18.184 17.200 18.167 /
\plot 17.200 18.167 17.173 18.150 /
\plot 17.173 18.150 17.141 18.133 /
\plot 17.141 18.133 17.109 18.117 /
\plot 17.109 18.117 17.075 18.100 /
\plot 17.075 18.100 17.041 18.085 /
\plot 17.041 18.085 17.005 18.068 /
\plot 17.005 18.068 16.967 18.051 /
\plot 16.967 18.051 16.925 18.034 /
\plot 16.925 18.034 16.880 18.017 /
\plot 16.880 18.017 16.836 18.000 /
\plot 16.836 18.000 16.789 17.983 /
\plot 16.789 17.983 16.749 17.968 /
\plot 16.749 17.968 16.715 17.956 /
\plot 16.715 17.956 16.690 17.947 /
\plot 16.690 17.947 16.675 17.941 /
\plot 16.675 17.941 16.669 17.939 /
%
\linethickness= 0.500pt
\setplotsymbol ({\thinlinefont .})
\setdashes < 0.1270cm>
\plot  7.620 24.130  6.985 24.130 /
%
\linethickness= 0.500pt
\setplotsymbol ({\thinlinefont .})
\plot  8.255 23.178  6.350 23.178 /
%
\linethickness= 0.500pt
\setplotsymbol ({\thinlinefont .})
\plot  8.890 22.225  5.715 22.225 /
%
\linethickness= 0.500pt
\setplotsymbol ({\thinlinefont .})
\plot  9.525 21.273  5.080 21.273 /
%
\linethickness= 0.500pt
\setplotsymbol ({\thinlinefont .})
\plot  8.890 20.320  5.715 20.320 /
%
\linethickness= 0.500pt
\setplotsymbol ({\thinlinefont .})
\plot  8.255 19.367  6.350 19.367 /
%
\linethickness= 0.500pt
\setplotsymbol ({\thinlinefont .})
\plot  7.620 18.415  6.985 18.415 /
%
\linethickness= 0.500pt
\setplotsymbol ({\thinlinefont .})
\setdots < 0.0953cm>
\plot 14.922 19.050 16.192 17.621 /
%
\linethickness= 0.500pt
\setplotsymbol ({\thinlinefont .})
\setsolid
\putrule from  5.715 22.225 to  5.713 22.225
\plot  5.713 22.225  5.704 22.221 /
\plot  5.704 22.221  5.685 22.214 /
\plot  5.685 22.214  5.654 22.204 /
\plot  5.654 22.204  5.609 22.187 /
\plot  5.609 22.187  5.560 22.170 /
\plot  5.560 22.170  5.512 22.151 /
\plot  5.512 22.151  5.465 22.134 /
\plot  5.465 22.134  5.423 22.115 /
\plot  5.423 22.115  5.385 22.100 /
\plot  5.385 22.100  5.351 22.085 /
\plot  5.351 22.085  5.319 22.070 /
\plot  5.319 22.070  5.292 22.056 /
\plot  5.292 22.056  5.266 22.039 /
\plot  5.266 22.039  5.239 22.024 /
\plot  5.239 22.024  5.213 22.007 /
\plot  5.213 22.007  5.188 21.988 /
\plot  5.188 21.988  5.163 21.969 /
\plot  5.163 21.969  5.137 21.948 /
\plot  5.137 21.948  5.110 21.927 /
\plot  5.110 21.927  5.084 21.903 /
\plot  5.084 21.903  5.059 21.880 /
\plot  5.059 21.880  5.033 21.859 /
\plot  5.033 21.859  5.008 21.836 /
\plot  5.008 21.836  4.985 21.812 /
\plot  4.985 21.812  4.964 21.791 /
\plot  4.964 21.791  4.942 21.770 /
\plot  4.942 21.770  4.921 21.749 /
\plot  4.921 21.749  4.900 21.728 /
\plot  4.900 21.728  4.879 21.706 /
\plot  4.879 21.706  4.858 21.685 /
\plot  4.858 21.685  4.837 21.662 /
\plot  4.837 21.662  4.813 21.637 /
\plot  4.813 21.637  4.792 21.611 /
\plot  4.792 21.611  4.769 21.586 /
\plot  4.769 21.586  4.750 21.560 /
\plot  4.750 21.560  4.729 21.533 /
\plot  4.729 21.533  4.712 21.507 /
\plot  4.712 21.507  4.695 21.482 /
\plot  4.695 21.482  4.680 21.457 /
\plot  4.680 21.457  4.667 21.431 /
\plot  4.667 21.431  4.657 21.404 /
\plot  4.657 21.404  4.646 21.378 /
\plot  4.646 21.378  4.638 21.351 /
\plot  4.638 21.351  4.631 21.321 /
\plot  4.631 21.321  4.625 21.292 /
\plot  4.625 21.292  4.621 21.260 /
\plot  4.621 21.260  4.616 21.226 /
\plot  4.616 21.226  4.614 21.192 /
\putrule from  4.614 21.192 to  4.614 21.160
\putrule from  4.614 21.160 to  4.614 21.126
\plot  4.614 21.126  4.616 21.095 /
\plot  4.616 21.095  4.619 21.065 /
\plot  4.619 21.065  4.621 21.035 /
\plot  4.621 21.035  4.625 21.008 /
\plot  4.625 21.008  4.631 20.980 /
\plot  4.631 20.980  4.638 20.951 /
\plot  4.638 20.951  4.644 20.921 /
\plot  4.644 20.921  4.655 20.892 /
\plot  4.655 20.892  4.665 20.862 /
\plot  4.665 20.862  4.678 20.832 /
\plot  4.678 20.832  4.691 20.803 /
\plot  4.691 20.803  4.705 20.775 /
\plot  4.705 20.775  4.720 20.750 /
\plot  4.720 20.750  4.737 20.724 /
\plot  4.737 20.724  4.754 20.703 /
\plot  4.754 20.703  4.771 20.682 /
\plot  4.771 20.682  4.790 20.663 /
\plot  4.790 20.663  4.807 20.646 /
\plot  4.807 20.646  4.828 20.629 /
\plot  4.828 20.629  4.849 20.612 /
\plot  4.849 20.612  4.875 20.595 /
\plot  4.875 20.595  4.900 20.578 /
\plot  4.900 20.578  4.928 20.561 /
\plot  4.928 20.561  4.957 20.546 /
\plot  4.957 20.546  4.987 20.532 /
\plot  4.987 20.532  5.016 20.517 /
\plot  5.016 20.517  5.046 20.504 /
\plot  5.046 20.504  5.076 20.491 /
\plot  5.076 20.491  5.108 20.479 /
\plot  5.108 20.479  5.133 20.468 /
\plot  5.133 20.468  5.160 20.458 /
\plot  5.160 20.458  5.190 20.447 /
\plot  5.190 20.447  5.220 20.436 /
\plot  5.220 20.436  5.251 20.426 /
\plot  5.251 20.426  5.283 20.413 /
\plot  5.283 20.413  5.315 20.403 /
\plot  5.315 20.403  5.347 20.392 /
\plot  5.347 20.392  5.376 20.384 /
\plot  5.376 20.384  5.406 20.373 /
\plot  5.406 20.373  5.433 20.367 /
\plot  5.433 20.367  5.459 20.358 /
\plot  5.459 20.358  5.482 20.352 /
\plot  5.482 20.352  5.503 20.345 /
\plot  5.503 20.345  5.531 20.339 /
\plot  5.531 20.339  5.556 20.335 /
\plot  5.556 20.335  5.582 20.331 /
\plot  5.582 20.331  5.609 20.326 /
\plot  5.609 20.326  5.637 20.324 /
\plot  5.637 20.324  5.664 20.322 /
\putrule from  5.664 20.322 to  5.690 20.322
\plot  5.690 20.322  5.707 20.320 /
\putrule from  5.707 20.320 to  5.713 20.320
\putrule from  5.713 20.320 to  5.715 20.320
%
\linethickness= 0.500pt
\setplotsymbol ({\thinlinefont .})
\plot  5.556 20.320  5.550 20.314 /
\plot  5.550 20.314  5.537 20.299 /
\plot  5.537 20.299  5.516 20.278 /
\plot  5.516 20.278  5.489 20.246 /
\plot  5.489 20.246  5.459 20.212 /
\plot  5.459 20.212  5.427 20.176 /
\plot  5.427 20.176  5.397 20.142 /
\plot  5.397 20.142  5.372 20.110 /
\plot  5.372 20.110  5.351 20.081 /
\plot  5.351 20.081  5.332 20.053 /
\plot  5.332 20.053  5.315 20.028 /
\plot  5.315 20.028  5.302 20.003 /
\plot  5.302 20.003  5.292 19.975 /
\plot  5.292 19.975  5.281 19.950 /
\plot  5.281 19.950  5.273 19.922 /
\plot  5.273 19.922  5.266 19.892 /
\plot  5.266 19.892  5.260 19.863 /
\plot  5.260 19.863  5.256 19.829 /
\plot  5.256 19.829  5.251 19.795 /
\plot  5.251 19.795  5.249 19.761 /
\putrule from  5.249 19.761 to  5.249 19.725
\putrule from  5.249 19.725 to  5.249 19.691
\plot  5.249 19.691  5.251 19.655 /
\plot  5.251 19.655  5.254 19.622 /
\plot  5.254 19.622  5.256 19.590 /
\plot  5.256 19.590  5.260 19.558 /
\plot  5.260 19.558  5.266 19.526 /
\plot  5.266 19.526  5.271 19.494 /
\plot  5.271 19.494  5.277 19.463 /
\plot  5.277 19.463  5.285 19.431 /
\plot  5.285 19.431  5.294 19.395 /
\plot  5.294 19.395  5.306 19.361 /
\plot  5.306 19.361  5.319 19.325 /
\plot  5.319 19.325  5.332 19.287 /
\plot  5.332 19.287  5.349 19.251 /
\plot  5.349 19.251  5.366 19.215 /
\plot  5.366 19.215  5.387 19.181 /
\plot  5.387 19.181  5.406 19.145 /
\plot  5.406 19.145  5.429 19.113 /
\plot  5.429 19.113  5.453 19.082 /
\plot  5.453 19.082  5.478 19.050 /
\plot  5.478 19.050  5.497 19.027 /
\plot  5.497 19.027  5.520 19.001 /
\plot  5.520 19.001  5.544 18.976 /
\plot  5.544 18.976  5.571 18.951 /
\plot  5.571 18.951  5.599 18.925 /
\plot  5.599 18.925  5.628 18.898 /
\plot  5.628 18.898  5.662 18.872 /
\plot  5.662 18.872  5.694 18.845 /
\plot  5.694 18.845  5.730 18.819 /
\plot  5.730 18.819  5.766 18.794 /
\plot  5.766 18.794  5.802 18.768 /
\plot  5.802 18.768  5.840 18.745 /
\plot  5.840 18.745  5.876 18.722 /
\plot  5.876 18.722  5.914 18.701 /
\plot  5.914 18.701  5.950 18.680 /
\plot  5.950 18.680  5.986 18.661 /
\plot  5.986 18.661  6.022 18.644 /
\plot  6.022 18.644  6.060 18.627 /
\plot  6.060 18.627  6.096 18.612 /
\plot  6.096 18.612  6.134 18.595 /
\plot  6.134 18.595  6.172 18.582 /
\plot  6.172 18.582  6.215 18.567 /
\plot  6.215 18.567  6.261 18.553 /
\plot  6.261 18.553  6.310 18.540 /
\plot  6.310 18.540  6.365 18.525 /
\plot  6.365 18.525  6.422 18.508 /
\plot  6.422 18.508  6.485 18.493 /
\plot  6.485 18.493  6.549 18.479 /
\plot  6.549 18.479  6.615 18.462 /
\plot  6.615 18.462  6.676 18.449 /
\plot  6.676 18.449  6.729 18.436 /
\plot  6.729 18.436  6.771 18.428 /
\plot  6.771 18.428  6.801 18.421 /
\plot  6.801 18.421  6.818 18.417 /
\plot  6.818 18.417  6.824 18.415 /
\putrule from  6.824 18.415 to  6.826 18.415
%
\linethickness= 0.500pt
\setplotsymbol ({\thinlinefont .})
\putrule from  6.826 24.130 to  6.824 24.130
\putrule from  6.824 24.130 to  6.818 24.130
\putrule from  6.818 24.130 to  6.799 24.130
\plot  6.799 24.130  6.767 24.128 /
\putrule from  6.767 24.128 to  6.723 24.128
\plot  6.723 24.128  6.670 24.126 /
\plot  6.670 24.126  6.612 24.124 /
\plot  6.612 24.124  6.553 24.122 /
\plot  6.553 24.122  6.498 24.117 /
\plot  6.498 24.117  6.447 24.113 /
\plot  6.447 24.113  6.401 24.109 /
\plot  6.401 24.109  6.356 24.105 /
\plot  6.356 24.105  6.318 24.100 /
\plot  6.318 24.100  6.284 24.094 /
\plot  6.284 24.094  6.251 24.086 /
\plot  6.251 24.086  6.219 24.077 /
\plot  6.219 24.077  6.187 24.067 /
\plot  6.187 24.067  6.155 24.056 /
\plot  6.155 24.056  6.124 24.043 /
\plot  6.124 24.043  6.090 24.031 /
\plot  6.090 24.031  6.058 24.016 /
\plot  6.058 24.016  6.026 23.999 /
\plot  6.026 23.999  5.992 23.980 /
\plot  5.992 23.980  5.961 23.961 /
\plot  5.961 23.961  5.929 23.939 /
\plot  5.929 23.939  5.899 23.918 /
\plot  5.899 23.918  5.870 23.897 /
\plot  5.870 23.897  5.840 23.876 /
\plot  5.840 23.876  5.814 23.853 /
\plot  5.814 23.853  5.789 23.832 /
\plot  5.789 23.832  5.764 23.808 /
\plot  5.764 23.808  5.743 23.785 /
\plot  5.743 23.785  5.719 23.764 /
\plot  5.719 23.764  5.696 23.738 /
\plot  5.696 23.738  5.675 23.713 /
\plot  5.675 23.713  5.652 23.688 /
\plot  5.652 23.688  5.630 23.658 /
\plot  5.630 23.658  5.609 23.628 /
\plot  5.609 23.628  5.588 23.599 /
\plot  5.588 23.599  5.567 23.567 /
\plot  5.567 23.567  5.548 23.535 /
\plot  5.548 23.535  5.529 23.501 /
\plot  5.529 23.501  5.512 23.470 /
\plot  5.512 23.470  5.497 23.438 /
\plot  5.497 23.438  5.484 23.404 /
\plot  5.484 23.404  5.472 23.372 /
\plot  5.472 23.372  5.461 23.340 /
\plot  5.461 23.340  5.450 23.309 /
\plot  5.450 23.309  5.442 23.277 /
\plot  5.442 23.277  5.433 23.245 /
\plot  5.433 23.245  5.427 23.209 /
\plot  5.427 23.209  5.423 23.175 /
\plot  5.423 23.175  5.417 23.137 /
\plot  5.417 23.137  5.412 23.099 /
\plot  5.412 23.099  5.410 23.061 /
\plot  5.410 23.061  5.408 23.023 /
\putrule from  5.408 23.023 to  5.408 22.985
\putrule from  5.408 22.985 to  5.408 22.947
\putrule from  5.408 22.947 to  5.408 22.911
\plot  5.408 22.911  5.410 22.875 /
\plot  5.410 22.875  5.412 22.843 /
\plot  5.412 22.843  5.417 22.811 /
\plot  5.417 22.811  5.419 22.782 /
\plot  5.419 22.782  5.425 22.754 /
\plot  5.425 22.754  5.431 22.718 /
\plot  5.431 22.718  5.438 22.684 /
\plot  5.438 22.684  5.446 22.650 /
\plot  5.446 22.650  5.457 22.617 /
\plot  5.457 22.617  5.469 22.585 /
\plot  5.469 22.585  5.480 22.553 /
\plot  5.480 22.553  5.493 22.523 /
\plot  5.493 22.523  5.505 22.496 /
\plot  5.505 22.496  5.518 22.473 /
\plot  5.518 22.473  5.531 22.449 /
\plot  5.531 22.449  5.544 22.428 /
\plot  5.544 22.428  5.556 22.409 /
\plot  5.556 22.409  5.571 22.388 /
\plot  5.571 22.388  5.586 22.369 /
\plot  5.586 22.369  5.603 22.348 /
\plot  5.603 22.348  5.624 22.324 /
\plot  5.624 22.324  5.645 22.299 /
\plot  5.645 22.299  5.671 22.274 /
\plot  5.671 22.274  5.692 22.250 /
\plot  5.692 22.250  5.707 22.233 /
\plot  5.707 22.233  5.713 22.227 /
\plot  5.713 22.227  5.715 22.225 /
%
\linethickness= 0.500pt
\setplotsymbol ({\thinlinefont .})
\putrule from  8.255 23.178 to  8.257 23.178
\plot  8.257 23.178  8.263 23.175 /
\plot  8.263 23.175  8.280 23.171 /
\plot  8.280 23.171  8.310 23.165 /
\plot  8.310 23.165  8.352 23.156 /
\plot  8.352 23.156  8.405 23.144 /
\plot  8.405 23.144  8.467 23.131 /
\plot  8.467 23.131  8.532 23.114 /
\plot  8.532 23.114  8.596 23.099 /
\plot  8.596 23.099  8.659 23.084 /
\plot  8.659 23.084  8.716 23.067 /
\plot  8.716 23.067  8.771 23.053 /
\plot  8.771 23.053  8.820 23.040 /
\plot  8.820 23.040  8.867 23.025 /
\plot  8.867 23.025  8.909 23.010 /
\plot  8.909 23.010  8.947 22.998 /
\plot  8.947 22.998  8.985 22.981 /
\plot  8.985 22.981  9.023 22.966 /
\plot  9.023 22.966  9.059 22.949 /
\plot  9.059 22.949  9.095 22.932 /
\plot  9.095 22.932  9.131 22.913 /
\plot  9.131 22.913  9.167 22.892 /
\plot  9.167 22.892  9.205 22.871 /
\plot  9.205 22.871  9.241 22.847 /
\plot  9.241 22.847  9.279 22.822 /
\plot  9.279 22.822  9.315 22.797 /
\plot  9.315 22.797  9.351 22.769 /
\plot  9.351 22.769  9.387 22.741 /
\plot  9.387 22.741  9.419 22.714 /
\plot  9.419 22.714  9.453 22.686 /
\plot  9.453 22.686  9.483 22.657 /
\plot  9.483 22.657  9.510 22.629 /
\plot  9.510 22.629  9.538 22.600 /
\plot  9.538 22.600  9.561 22.572 /
\plot  9.561 22.572  9.584 22.545 /
\plot  9.584 22.545  9.605 22.515 /
\plot  9.605 22.515  9.627 22.483 /
\plot  9.627 22.483  9.648 22.451 /
\plot  9.648 22.451  9.667 22.416 /
\plot  9.667 22.416  9.684 22.380 /
\plot  9.684 22.380  9.701 22.341 /
\plot  9.701 22.341  9.716 22.303 /
\plot  9.716 22.303  9.730 22.263 /
\plot  9.730 22.263  9.741 22.221 /
\plot  9.741 22.221  9.751 22.181 /
\plot  9.751 22.181  9.758 22.138 /
\plot  9.758 22.138  9.764 22.098 /
\plot  9.764 22.098  9.768 22.058 /
\putrule from  9.768 22.058 to  9.768 22.020
\putrule from  9.768 22.020 to  9.768 21.982
\plot  9.768 21.982  9.766 21.943 /
\plot  9.766 21.943  9.764 21.907 /
\plot  9.764 21.907  9.758 21.872 /
\plot  9.758 21.872  9.749 21.833 /
\plot  9.749 21.833  9.741 21.795 /
\plot  9.741 21.795  9.728 21.753 /
\plot  9.728 21.753  9.713 21.709 /
\plot  9.713 21.709  9.694 21.660 /
\plot  9.694 21.660  9.673 21.607 /
\plot  9.673 21.607  9.650 21.550 /
\plot  9.650 21.550  9.624 21.491 /
\plot  9.624 21.491  9.597 21.431 /
\plot  9.597 21.431  9.574 21.378 /
\plot  9.574 21.378  9.553 21.334 /
\plot  9.553 21.334  9.538 21.300 /
\plot  9.538 21.300  9.529 21.281 /
\plot  9.529 21.281  9.525 21.275 /
\putrule from  9.525 21.275 to  9.525 21.273
%
\linethickness= 0.500pt
\setplotsymbol ({\thinlinefont .})
\plot  9.525 21.273  9.531 21.268 /
\plot  9.531 21.268  9.546 21.258 /
\plot  9.546 21.258  9.569 21.241 /
\plot  9.569 21.241  9.601 21.220 /
\plot  9.601 21.220  9.639 21.192 /
\plot  9.639 21.192  9.682 21.160 /
\plot  9.682 21.160  9.724 21.129 /
\plot  9.724 21.129  9.764 21.099 /
\plot  9.764 21.099  9.800 21.071 /
\plot  9.800 21.071  9.832 21.044 /
\plot  9.832 21.044  9.862 21.018 /
\plot  9.862 21.018  9.887 20.995 /
\plot  9.887 20.995  9.908 20.974 /
\plot  9.908 20.974  9.929 20.951 /
\plot  9.929 20.951  9.948 20.927 /
\plot  9.948 20.927  9.970 20.902 /
\plot  9.970 20.902  9.989 20.875 /
\plot  9.989 20.875 10.005 20.845 /
\plot 10.005 20.845 10.022 20.813 /
\plot 10.022 20.813 10.039 20.781 /
\plot 10.039 20.781 10.056 20.745 /
\plot 10.056 20.745 10.069 20.709 /
\plot 10.069 20.709 10.082 20.673 /
\plot 10.082 20.673 10.094 20.635 /
\plot 10.094 20.635 10.105 20.597 /
\plot 10.105 20.597 10.113 20.561 /
\plot 10.113 20.561 10.122 20.525 /
\plot 10.122 20.525 10.128 20.489 /
\plot 10.128 20.489 10.135 20.451 /
\plot 10.135 20.451 10.137 20.424 /
\plot 10.137 20.424 10.141 20.394 /
\plot 10.141 20.394 10.143 20.364 /
\plot 10.143 20.364 10.145 20.333 /
\putrule from 10.145 20.333 to 10.145 20.301
\plot 10.145 20.301 10.143 20.267 /
\plot 10.143 20.267 10.139 20.231 /
\plot 10.139 20.231 10.135 20.197 /
\plot 10.135 20.197 10.126 20.161 /
\plot 10.126 20.161 10.118 20.125 /
\plot 10.118 20.125 10.105 20.091 /
\plot 10.105 20.091 10.090 20.055 /
\plot 10.090 20.055 10.073 20.022 /
\plot 10.073 20.022 10.054 19.990 /
\plot 10.054 19.990 10.031 19.958 /
\plot 10.031 19.958 10.008 19.928 /
\plot 10.008 19.928  9.978 19.899 /
\plot  9.978 19.899  9.948 19.869 /
\plot  9.948 19.869  9.921 19.848 /
\plot  9.921 19.848  9.891 19.825 /
\plot  9.891 19.825  9.857 19.801 /
\plot  9.857 19.801  9.823 19.778 /
\plot  9.823 19.778  9.785 19.755 /
\plot  9.785 19.755  9.743 19.732 /
\plot  9.743 19.732  9.701 19.710 /
\plot  9.701 19.710  9.656 19.687 /
\plot  9.656 19.687  9.610 19.664 /
\plot  9.610 19.664  9.561 19.643 /
\plot  9.561 19.643  9.510 19.619 /
\plot  9.510 19.619  9.459 19.600 /
\plot  9.459 19.600  9.409 19.579 /
\plot  9.409 19.579  9.358 19.560 /
\plot  9.358 19.560  9.309 19.541 /
\plot  9.309 19.541  9.260 19.524 /
\plot  9.260 19.524  9.212 19.509 /
\plot  9.212 19.509  9.165 19.494 /
\plot  9.165 19.494  9.121 19.482 /
\plot  9.121 19.482  9.078 19.469 /
\plot  9.078 19.469  9.036 19.456 /
\plot  9.036 19.456  8.996 19.446 /
\plot  8.996 19.446  8.947 19.435 /
\plot  8.947 19.435  8.903 19.425 /
\plot  8.903 19.425  8.856 19.416 /
\plot  8.856 19.416  8.810 19.408 /
\plot  8.810 19.408  8.763 19.401 /
\plot  8.763 19.401  8.712 19.395 /
\plot  8.712 19.395  8.661 19.391 /
\plot  8.661 19.391  8.606 19.384 /
\plot  8.606 19.384  8.549 19.380 /
\plot  8.549 19.380  8.492 19.378 /
\plot  8.492 19.378  8.435 19.374 /
\plot  8.435 19.374  8.384 19.372 /
\plot  8.384 19.372  8.338 19.370 /
\putrule from  8.338 19.370 to  8.302 19.370
\plot  8.302 19.370  8.276 19.367 /
\putrule from  8.276 19.367 to  8.261 19.367
\putrule from  8.261 19.367 to  8.255 19.367
\linethickness=0pt
\putrectangle corners at  4.589 24.632 and 18.864 17.596
\endpicture}

\end{center}
\vspace{.5cm}

An alternate construction of a space of complete skew forms has been
given by Bertram \cite{bert}.  Despite his different approach using
Pfaffians, he obtains what appear to be the same spaces as those
discussed here.

\bit{Appendix: symplectic quotients and broken Morse flows
  \label{symplectic}} 

At least heuristically, the Chow quotient of the Grassmannian may be
understood in terms of broken Morse flows.  Without proving anything,
this final section will explain that analogy.

Let $X$ be a smooth compact manifold with a symplectic form $\omega$,
and a Hamiltonian action of the Lie group $\U{1}$.  Here {\em
Hamiltonian} means that the action preserves $\omega$, and admits a
{\em moment map}, namely a primitive function for the infinitesimal
action in the following sense.  The infinitesimal action of $\U{1}$
determines a vector field $\xi \in \Ga(TX)$.  A moment map is a smooth
map $\mu: X \to \R$ such that $d \mu = \omega(\xi,\cdot)$ as
$1$-forms.  Assume for simplicity that the action is {\em quasi-free},
meaning that it has no nontrivial finite stabilizers.  If $t$ is a
regular value of $\mu$, then the {\em symplectic quotient} or {\em
  reduced phase space} $\mu^{-1}(t) \mod \U{1}$ is a manifold, which
turns out to have a natural symplectic structure.  (It first arose in
classical mechanics, where $X$ was a phase space, the group action
reflected some symmetry of the system, and the moment map was the
associated conserved quantity.)

Several authors \cite{gs,hu} have studied the relationship between the
symplectic quotients for different values of $t$.  A key point is that
$\mu$ is a Morse function in the sense of Bott \cite{bott}.  Indeed,
the critical manifolds are precisely the fixed manifolds of the
$\U{1}$-action.  The different level sets are therefore related by a
sequence of surgeries: submanifolds diffeomorphic to sphere bundles
over the critical manifolds are removed, and replaced by sphere
bundles of different dimension.  The $\U{1}$-action preserves these
bundles, and indeed acts in the standard way on the odd-dimensional
spheres which are the fibers.  Dividing by this action produces
bundles with fiber a complex projective space, so the quotients are
related by excising one projective bundle and replacing it with
another.

Those who prefer a quotient to be a canonical object may be unhappy
that the symplectic quotient depends on the choice of $t$.  They will
prefer the following alternative.  Let $F$ be the space of Morse flow
lines from the absolute maximum of $\mu$ to its absolute minimum.
This is a manifold which can clearly be identified with a dense open
subset of every level set, simply by intersecting the level set with
the flow.  It is well-known in Morse theory \cite{ab} that $F$ can be
compactified by adding the so-called {\em broken flows}.  These are
unions of flow lines, each beginning at the same critical point where
the previous one ends.  The picture is exactly like the one shown in
\S\ref{morphism}.  

The space $\overline{F}$ of broken flows is not generally a manifold,
but only a manifold with corners.  The $\U{1}$-action on
$\overline{F}$ is free, so the quotient is again a manifold with
corners.  But each stratum is foliated by tori coming from the $\U{1}$
actions on the components of the broken flow: the smaller the stratum,
the larger the torus.  If these tori are all collapsed, then the
resulting space $\X$ is a symplectic manifold.  It parametrizes broken
flows modulo an equivalence in which each component of the flow
rotates separately.  The proof that $\X$ is smooth is not easy to find
in the literature, but it resembles the line of argument in Lerman's
work on symplectic cuts \cite{l}.

Another way to regard this space $\X$ is to use the Morse flow to
induce maps from the quotient at each regular value to the quotients
at the two adjacent critical values.  This makes the symplectic
quotients into an inverse system, and $\X$ is the inverse limit.

In the body of the paper we were especially interested in an action of
$\cx$ on the Grassmannian $\Gr_u(V\oplus W)$.  Although we won't go
into the details of the proof, it follows from the well-known
identification of symplectic quotients with Mumford quotients
\cite[\S8.2]{mfk} that in such a case $\X$ is symplectomorphic to the
inverse limit of Mumford quotients referred to in \re{gr}{c} of the
Main Theorem.  The $\U{1}$ orbits of Morse flow lines get replaced by
$\cx$-orbit closures, so equivalence classes of broken flows get
replaced by nodal curves of genus 0: exactly those described in
Proposition \ref{nodal}.  The symplectic point of view probably will
not imply any new results here.  Nevertheless, thinking of points in
the limit of Mumford quotients, or in the Chow quotient, as broken
flows is a valuable aid to the intuition.

\end{document}